\documentclass[11pt,letterpaper,reqno]{amsart}
\usepackage{amsmath,amsthm,amsfonts,amssymb,mathtools,algpseudocode}
\usepackage{comment,bookmark,bm,color}
\usepackage{comment}
\hypersetup{pdfstartview={FitH}}

\newtheorem{Theorem}{{\bf Theorem}}[section]

\newtheorem{Corollary}[Theorem]{{\bf Corollary}}
\newtheorem{Algorithm}[Theorem]{{\bf Algorithm}}
\newtheorem{Proposition}[Theorem]{{\bf Proposition}}
\newtheorem{Definition}[Theorem]{{\bf Definition}}
\newtheorem{Lemma}[Theorem]{{\bf Lemma}}

\numberwithin{equation}{section}

\newcommand{\ve}{\text{vec}}

\newcommand{\diag}{\text{diag}}

\newcommand{\Rl}{\mbox{\Large $\mathcal{R}$}}
\newcommand{\Ol}{\mbox{\Large $\mathcal{O}$}}
\newcommand{\Jl}{\mbox{\Large $\mathcal{J}$}}
\newcommand\restrict[1]{\raisebox{-.5ex}{$|$}_{#1}}

\begin{document}

\title[On the convergence of complex Jacobi methods]{On the convergence of complex Jacobi methods}
\author{Vjeran Hari}\thanks{Vjeran Hari, Department of Mathematics, Faculty of Science, University of Zagreb, Bijeni\v{c}ka 30, 10000 Zagreb, Croatia, \texttt{hari@math.hr}}
\author{Erna Begovi\'{c}~Kova\v{c}}\thanks{Erna Begovi\'{c} Kova\v{c}, Faculty of Chemical Engineering and Technology, University of Zagreb, Maruli\'{c}ev trg 19, 10000 Zagreb, Croatia, \texttt{ebegovic@fkit.hr}}

\thanks{This work has been fully supported by Croatian Science Foundation under the project 3670.}
\date{\today}


\begin{abstract}
In this paper we prove the global convergence of the complex Jacobi method for Hermitian matrices for a large class of generalized serial pivot strategies. For a given Hermitian matrix $A$ of order $n$ we find a constant $\gamma<1$ depending on $n$, such that $S(A')\leq\gamma{S(A)}$, where $A'$ is obtained from $A$ by applying one or more cycles of the Jacobi method and $S(\cdot)$ stands for the off-norm. Using the theory of complex Jacobi operators, the result is generalized so it can be used for proving convergence of more general Jacobi-type processes. In particular, we use it to prove the global convergence of Cholesky-Jacobi method for solving the positive definite generalized eigenvalue problem.
\end{abstract}

\maketitle

\section{Introduction}

Let $A$ be a Hermitian matrix of order $n$. Complex Jacobi method is the iterative process
\begin{equation}\label{hjacobiagm}
A^{(k+1)}=U_k^*A^{(k)}U_k, \quad k\geq0, \quad A^{(0)}=A,
\end{equation}
where $U_k=R(i_k,j_k,\phi_k,\alpha_k)$ are complex plane rotations of the form
{\footnotesize
\begin{equation}\label{R}
R(i_k,j_k,\phi_k,\alpha_k)=\left[
    \begin{array}{ccccc}
       I &  &  &  & \\
       & \cos\phi_k &  & -e^{\imath\alpha_k}\sin\phi_k & \\
       &  & I &  & \\
       & e^{-\imath\alpha_k}\sin\phi_k &  & \cos\phi_k & \\
       &  &  &  & I \\
    \end{array}
  \right]
  \begin{array}{l}
     \\
     i_k \\
     \\
     j_k \\
     \\
     \end{array},
\end{equation}}
with $\imath=\sqrt{-1}$. The angles $\phi_k$ and $\alpha_k$ are chosen to annihilate the pivot element $a_{i_kj_k}^{(k)}$ of $A^{(k)}$. This choice of angles maximally reduces $S(A^{(k+1)})$ for some pivot pair $(i_k,j_k)$. Here
\[
S(X)=\|X-\diag(X)\|_F
\]
stands for the off-norm of matrix $X$, which is the Frobenius norm of the off-diagonal part of $X$.

The real symmetric Jacobi method is known for its efficiency~\cite{drm+ves-04a,drm+ves-04b}, high relative accuracy~\cite{dem+ves-92} and convergence properties~\cite{fer-89,har-91,har+beg-17}. Among all this, it is very suitable for parallel computation~\cite{sam-71,luk+par-89} which enhanced its usage in modern CPU and GPU computation \cite{nov+sin-15}.
Most of these properties, especially the parallelism, are also inherent in the complex Jacobi method.
However, the global convergence theory of the symmetric Jacobi method, which uses Jacobi annihilators and operators \cite{zim-65,har+beg-17}, cannot be straightforwardly applied to the complex Hermitian Jacobi method (cf. \cite[Section~2]{har-86}).

In this paper we study the global convergence of the complex Jacobi method under the large class of generalized serial strategies from~\cite{har+beg-17}. To this end we define a more general iterative process described by complex Jacobi annihilators and operators that generalize one step and one cycle (sweep) of the Jacobi method. Compared to the real Jacobi annihilators and operators from~\cite{hen+zim-68,zim-65}, here we have to do some adjustments, since we work with complex Hermitian matrices instead of symmetric ones. Our convergence result takes the form
\[
S(A')\leq\gamma{S(A)}, \quad 0\leq\gamma<1,
\]
where $A'$ is obtained from $A$ by applying one or more cycles of the Jacobi method. The constant $\gamma$ does not depend on a particular matrix $A$, only on its size. This paper extends the known Forsythe-Henrici condition from~\cite{for+hen-60} to the generalized serial complex Jacobi methods.

In the second part of the paper we use the results obtained for the complex Jacobi annihilators and operators to prove the global convergence of a Jacobi method for the positive definite generalized eigenvalue problem $Ax=\lambda Bx$, where $A$ and $B$ are Hermitian and $B$ is positive definite. This is a rather new method, introduced in~\cite{har-AIP} (the real method is analyzed in \cite{har-HZ}) and is known for its outstanding relative accuracy.

The paper is divided in four sections. In Section~\ref{sec:complJac} we define and analyze the complex Jacobi method for a Hermitian matrix $A$ and introduce the corresponding Jacobi annihilators and operators.
In Section~\ref{sec:strategies} we study the global convergence of that method under the large class of generalized serial  pivot strategies. In Section~\ref{sec:pgep} we prove the global convergence of the Cholesky-Jacobi method for the positive definite generalized eigenvalue problem.

\section{Complex Jacobi method}\label{sec:complJac}

We briefly describe one step of the complex Jacobi method for Hermitian matrices. Then we define the Jacobi annihilators and operators and prove some auxiliary results which are later used in the global convergence analysis.

In the $k$th step, the angles $\alpha_k$ and $\phi_k$ from \eqref{R} are chosen to annihilate the elements $a_{i_k j_k}$ and $a_{j_k i_k}$. To simplify notation, let us fix $k$ and use $A=(a_{rt})$ and $A'=(a_{rt}')$ for the current and updated iteration matrix, $(i,j)$ for the pivot pair $(i_k,j_k)$, and $\alpha$, $\phi$ for the angles $\alpha_k$, $\phi_k$, respectively. Similarly, the complex rotation $U_k$ is denoted by $U$.
The $2\times 2$ pivot submatrices $\hat{A}$, $\hat{U}$ are obtained on the intersection of pivot rows and columns of $A$ and $U$, respectively.

The values of $\alpha$ and $\phi$ are determined from the relation
\begin{equation}\label{submatrix}
\null\hspace{-2ex} \left[
    \begin{array}{cc}
      a_{ii}' & 0 \\
      0 & a_{jj}' \\
    \end{array}  \right]
    =
    \hat{U}^*\hat{A}\hat{U}
    =
    \left[ \!
    \begin{array}{cc}
      \cos\phi & -e^{\imath\alpha}\sin\phi \\
      e^{-\imath\alpha}\sin\phi & \cos\phi \\
    \end{array}
  \right]^*\left[
    \begin{array}{cc}
      a_{ii} & a_{ij} \\
      a_{ji} & a_{jj} \\
    \end{array}
  \right]\left[
    \begin{array}{cc}
      \cos\phi & -e^{\imath\alpha}\sin\phi \\
      e^{-\imath\alpha}\sin\phi & \cos\phi \\
    \end{array}
  \right  ]\! .
\end{equation}
After multiplying the equation~\eqref{submatrix} by $\hat{U}$ from the left side, because the diagonal elements on the left- and right-hand side are equal, we obtain
\begin{align*}
a_{ii}'= & a_{ii}+a_{ij}e^{-\imath\alpha}\tan\phi,\\
a_{jj}'= & a_{jj}-a_{ji}e^{\imath\alpha}\tan\phi.
\end{align*}
Since the trace is invariant to unitary transformations, we have
$$a_{ii}+a_{jj} = a_{ii}'+a_{jj}' = a_{ii}+a_{ij}e^{-\imath\alpha}\tan\phi+a_{jj}-a_{ji}e^{\imath\alpha}\tan\phi,$$
that is
$$a_{ij}e^{-\imath\alpha}=a_{ji}e^{\imath\alpha}.$$
Matrix $\hat{A}$ is Hermitian, hence we can write $a_{ij}=|a_{ij}|e^{\imath\text{arg}(a_{ij})}$, $a_{ji}=|a_{ij}|e^{-\imath\text{arg}(a_{ij})}$. We obtain
$|a_{ij}|e^{ i \alpha_{ij}}e^{- i \alpha} = |a_{ij}|e^{- i \alpha_{ij}}e^{ i \alpha}$, which implies
\begin{equation}\label{alpha}
\alpha=\text{arg}(a_{ij}).
\end{equation}
Further on, from~\eqref{submatrix} we have
$$0=a_{ij}' = \left[
    \begin{array}{cc}
      \cos\phi & e^{\imath\alpha}\sin\phi \\
    \end{array}
  \right]\left[
    \begin{array}{cc}
      a_{ii} & a_{ij} \\
      a_{ji} & a_{jj} \\
    \end{array}
  \right]\left[
    \begin{array}{c}
      -e^{\imath\alpha}\sin\phi \\
      \cos\phi \\
    \end{array}
  \right].$$
Thus, using~\eqref{alpha} we get the known formula
$$\tan(2\phi)=\frac{2|a_{ij}|}{a_{ii}-a_{jj}}.$$
A choice for the angle interval is $[-\pi/4,\pi /4]$. To compute $\cos\phi$ and $\sin\phi$
one can choose stable formulas
\begin{equation}\label{hjacobit}
t = \frac{2|a_{ij}|\text{sign}(a_{ii}-a_{jj})}{|a_{ii}-a_{jj}|+\sqrt{|a_{ii}-a_{jj}|^2+4|a_{ij}|^2}},
\quad \cos\phi = \frac{1}{\sqrt{1+t^2}},\quad \sin\phi = \frac{t}{\sqrt{1+t^2}},
\end{equation}
where $t$ denotes $\tan\phi$.
One iteration of the complex Jacobi method is presented in Algorithm~\ref{agm:hJacobi}. The algorithm uses only the upper-triangle of the matrix.

\begin{Algorithm}\label{agm:hJacobi}
\hrule\vspace{1ex}
\emph{One iteration of the complex Jacobi algorithm}
\vspace{0.5ex}\hrule
\begin{algorithmic}
\If {$a_{ij}\neq0$}
\State $t = 2|a_{ij}|\text{sign}(a_{ii}-a_{jj})/(|a_{ii}-a_{jj}|+\sqrt{|a_{ii}-a_{jj}|^2+4|a_{ij}|^2})$;
\State $c=1/\sqrt{1+t^2}$; \quad $s=t/\sqrt{1+t^2}$;  \quad $s^+=s\,a_{ij}/|a_{ij}|$; \quad $s^-=s\,|a_{ij}|/a_{ij};$
\State $a_{ii}=a_{ii}+t|a_{ij}|$; \quad $a_{jj}=a_{jj}-t|a_{ij}|$; \quad $a_{ij}=0$; \quad $a_{ji}=0;$
\For {$r=1,\ldots,i-1$}
\State $x=ca_{ri}+s^- a_{rj}; \quad a_{rj}=-s^+a_{ri}+ca_{rj}; \quad a_{ri}=x;$
\EndFor
\For {$r=i+1,\ldots,j-1$}
\State $x=ca_{ir}+s^+\bar{a}_{rj}; \quad a_{rj}=-s^+\bar{a}_{ir}+ca_{rj}; \quad a_{ir}=x;$
\EndFor
\For {$r=j+1,\ldots,n$}
\State $x=ca_{ir}+s^+a_{jr}; \quad a_{jr}=-s^-a_{ir}+ca_{jr}; \quad a_{ir}=x;$
\EndFor
\EndIf
\end{algorithmic}
\hrule
\end{Algorithm}
\bigskip

A Jacobi method is \emph{convergent} on $A$ if the sequence of matrices $(A^{(k)}, \ k\geq0)$ obtained by the iterative process~\eqref{hjacobiagm} converges to some diagonal matrix. The method is \emph{globally convergent} if it is convergent on every Hermitian matrix $A$.

A \emph{pivot strategy} can be identified with a function
$I:\mathbb{N}_0\rightarrow \mathcal{P}_n,$
where $\mathbb{N}_0=\{0,1,2,3,\ldots\}$ and $\mathcal{P}_n=\{(r, s) \,|\, 1\leq r<s\leq n\}$.
If $I$ is a periodic function with period $T=N\equiv\frac{n(n-1)}{2}$ and $\{I(0),I(1),\ldots,I(T-1)\}=\mathcal{P}_n$, then $I$ is a
\emph{cyclic} pivot strategy.
By $\mathcal{\Ol}(\mathcal{S})$ we denote the \emph{set of all finite sequences} made of the elements of
$\mathcal{S}$, $\mathcal{S}\subseteq \mathcal{P}_n$, assuming that each pair from $\mathcal{S}$ appears at least once in each sequence.
If $I$ is a cyclic strategy then $\mathcal{O}_I$ stands for the sequence $I(0),I(1),\ldots,I(N-1)\in\mathcal{\Ol}(\mathcal{P}_n)$ generated by the first $N$ steps (i.e., by the first \emph{sweep}) of the method. Obviously, $\mathcal{O}_I$ is an ordering of $\mathcal{P}_n$ and it is sometimes called \emph{pivot ordering}.
Conversely, if $\mathcal{O}\in\mathcal{\Ol}(\mathcal{P}_n)$ is an ordering of $\mathcal{P}_n$,  $\mathcal{O}=(i_0,j_0),(i_1,j_1),\ldots,(i_{N-1},j_{N-1})$ then the cyclic strategy $I_{\mathcal{O}}$ is defined by $I_{\mathcal{O}}(k)=(i_{\tau(k)},j_{\tau(k)})$, where $\tau(k)$ is the unique integer satisfying $0\leq \tau(k)\leq N-1$ and $k\equiv\tau(k)(\!\!\! \mod \ N)$, $k\geq0$.

We define several equivalence relations on the set $\mathcal{\Ol}(\mathcal{S})$, $\mathcal{S}\subseteq\mathcal{P}_n$.

An \emph{admissible transposition} in $\mathcal{O}\in\mathcal{\Ol}(\mathcal{S})$ is any transposition of two adjacent terms in $\mathcal{O}$,
\[
(i_r,j_r),(i_{r+1},j_{r+1})\rightarrow(i_{r+1},j_{r+1}),(i_r,j_r),
\]
provided that the sets $\{i_r,j_r\}$ and $\{i_{r+1},j_{r+1}\}$ are disjoint. We also say that the pairs
$(i_r,j_r)$ and $(i_{r+1},j_{r+1})$ are disjoint or that they \emph{commute}. The length of $\mathcal{O}$
is the number of pairs contained in it. The following definition is taken from \cite{har+beg-17}.

\begin{Definition}
Two sequences $\mathcal{O},\mathcal{O}'\in\mathcal{\Ol}(\mathcal{S})$, $\mathcal{S}\subseteq \mathcal{P}_n$, are said to be
\begin{itemize}
\item[(i)] \emph{equivalent} (we write $\mathcal{O}\sim\mathcal{O}'$) if one can be obtained from the other by a finite set of admissible transpositions
\item[(ii)] \emph{shift-equivalent} ($\mathcal{O}\stackrel{\mathsf{s}}{\sim}\mathcal{O}'$) if $\mathcal{O}=[\mathcal{O}_1,\mathcal{O}_2]$ and $\mathcal{O}'=[\mathcal{O}_2,\mathcal{O}_1]$, where $[ \ , \  ]$ stands for concatenation; the length of $\mathcal{O}_1$ is called shift length
\item[(iii)] \emph{weak equivalent} ($\mathcal{O}\stackrel{\mathsf{w}}{\sim}\mathcal{O}'$) if there exist $\mathcal{O}_i\in\mathcal{\Ol}(\mathcal{S})$, $0\leq i\leq r$, such that every two adjacent terms in the sequence $\mathcal{O}=\mathcal{O}_0,\mathcal{O}_1,\ldots,\mathcal{O}_r=\mathcal{O}'$ are equivalent or shift-equivalent
\item[(iv)] \emph{permutation equivalent} ($\mathcal{O}\stackrel{\mathsf{p}}{\sim}\mathcal{O}'$ or $\mathcal{O}'=\mathcal{O}(\mathsf{q})$) if there is a permutation $\mathsf{q}$ of the set $\mathcal{S}$ such that
    $\mathcal{O}'=(\mathsf{q}(i_0),\mathsf{q}(j_0)),(\mathsf{q}(i_1),\mathsf{q}(j_1)),\ldots,
    (\mathsf{q}(i_r),\mathsf{q}(j_r))$
\item[(v)] \emph{reverse} ($\mathcal{O}'=\mathcal{O}^{\leftarrow}$) if $\mathcal{O}=(i_0,j_0),(i_1,j_1),\ldots,(i_{r},j_{r})$ and $\mathcal{O}'=(i_{r},j_{r}),\ldots,(i_1,j_1),(i_0,j_0)$.
\end{itemize}
\end{Definition}

Two pivot strategies $I_{\mathcal{O}}$ and $I_{\mathcal{O}'}$ are equivalent (resp.\@ shift-equivalent, weak equivalent, permutation equivalent, reverse) if the corresponding sequences $\mathcal{O}$ and $\mathcal{O}'$ are equivalent (resp.\@ shift-equivalent, weak equivalent, permutation equivalent, reverse).
If $\mathcal{O}\stackrel{\mathsf{w}}{\sim}\mathcal{O}'$, then there is a sequence $\mathcal{O}=\mathcal{O}_0,\mathcal{O}_1,\ldots,\mathcal{O}_r=\mathcal{O}'$
such that
\begin{equation}\label{canonical}
\mathcal{O}_0\sim
\mathcal{O}_1\stackrel{\mathsf{s}}{\sim}\mathcal{O}_2\sim\mathcal{O}_3\stackrel{\mathsf{s}}{\sim}\mathcal{O}_4\ldots \qquad \text{or} \qquad \mathcal{O}_0\stackrel{\mathsf{s}}{\sim}
\mathcal{O}_1\sim\mathcal{O}_2\stackrel{\mathsf{s}}{\sim}\mathcal{O}_3\sim\mathcal{O}_4\ldots.
\end{equation}
We say that the chain from~\eqref{canonical} connecting $\mathcal{O}$ and $\mathcal{O}'$ is in the \emph{canonical form}.

\subsection{Jacobi annihilators and operators}

Jacobi annihilators and operators are the tools for analyzing the convergence of cyclic and quasi-cyclic Jacobi methods. They appear when the iteration matrices are represented by vectors.
They were introduced in~\cite{zim-65,hen+zim-68} and later were generalized to work with complex and block  methods (see~\cite{har-86,har-09,har-15,har+beg-17}).
As we will see later, Jacobi annihilators and operators are crucial in proving the convergence of the method for some cyclic pivot strategies.

For $2\leq s,t\leq n$ we define vectors (cf. \cite{har-86})
\[
c_t = \left[
          \begin{array}{c}
            a_{1t} \\
            a_{2t} \\
            \vdots \\
            a_{t-1,t} \\
          \end{array}
        \right], \qquad r_s=[a_{s1},\ldots,a_{s,s-1}]
        \]
and function $\ve:\mathbb{C}^{n\times n}\rightarrow \mathbb{C}^{2N}$,
\[
\ve(A)=[c_2^T,c_3^T,\ldots,c_n^T,r_2,r_3,\ldots,r_n]^T.
\]
It is easy to see that $\ve $ is a linear operator and a surjection.
If $A$ is Hermitian, we have
\begin{equation}\label{herm_A}
a=\ve(A)=\left[
             \begin{array}{c}
               v \\
               \bar{v} \\
             \end{array}
           \right],
\end{equation}
where $v=[c_2^T,c_3^T,\ldots,c_n^T]^T$. In that case, $\emph{\ve}(A)$ is determined by the strict upper-triangular part of $A$.

Let $\nu_{ij}:\mathbb{C}^{n\times n}\rightarrow\mathbb{C}^{n\times n}$ be the linear operator that sets the elements on pivot positions $(i,j)$ and $(j,i)$ in matrix $A$ to zero.

\begin{Definition}\label{def: JacAnn}
Let $U=R(i,j,\phi,\alpha)$ be a complex rotation. The matrix $\mathcal{R}_{ij}(U)$ defined by
\begin{equation}\label{def:jan}
\mathcal{R}_{ij}(U)\,\emph{\ve}(A) = \emph{\ve}(\nu_{ij}(U^*AU)), \quad A\in\mathbb{C}^{n\times n},
\end{equation}
is called Jacobi annihilator.
The class of the Jacobi annihilators $\Rl_{ij}^{\nu }$, $\nu\in [0,1]$ is a set
\[
\Rl_{ij}^{\nu } = \Big{\{} \mathcal{R}_{ij}(U) \ \Big{|} \ U=R(i,j,\phi,\alpha), \ 0\leq\phi,\alpha\leq 2\pi,\
|cos\phi|\geq \nu \Big{\}}.
\]
If $\nu =0$,  $\Rl_{ij}$ is used instead of $\Rl_{ij}^{0}$.
\end{Definition}

As it is suggested by the notation, the Jacobi annihilator depends on $U$ and the pair $(i,j)$. Since $\emph{\ve}$ is surjection, $\emph{\ve}(A)$ passes through all elements of $\mathbb{C}^{2N}$, so $\mathcal{R}_{ij}(U)$ is well defined. More on $\ve$ can be found in \cite[Section~2]{har-09}, \cite[Section~4]{har-15}.
Since $\ve$ is not invertible, we shall also use $\ve_0 =\ve\restrict{\mathbb{C}_0^{n\times n}}$ where
$\mathbb{C}_0^{n\times n}$ is the set of all complex matrices of dimension $n$ with zero diagonal elements. Therefore, if we are given $a\in \mathbb{C}^{2N}$, $A=\ve_0^{-1}(a)$ is a complex matrix of dimension $n$ with zero diagonal elements and the first (second) $N$ elements of $a$ determine the upper-(lower-)triangle of $A$ in such a way that $\ve(\ve_0(a))=a$.

The parameter $\nu$ is linked to the global convergence of Jacobi methods. The known Forsythe-Henrici condition \cite{for+hen-60} for the convergence of the serial Jacobi methods requires $\nu > 0$. In reality, the complex Jacobi method for Hermitian matrices will always satisfy $\cos\phi_k \geq \sqrt{2}/2$, hence for that method we have $\nu=\sqrt{2}/2$. However, if the initial Hermitian matrix is nearly diagonal, $\nu$ can be larger than $\sqrt{2}/2$. As we are going to see in Section~4, an appropriate $\nu$ will play important role in the global convergence proof of the Cholesky-Jacobi method for the positive definite generalized eigenvalue problem.

In the both relations~\eqref{hjacobiagm} and~\eqref{def:jan} the plane rotation is applied to the underlaying matrix through similarity transformation.
The difference is that the rotation angles in~\eqref{hjacobiagm} are chosen to annihilate the pivot element, while in~\eqref{def:jan} one can take arbitrary rotation angles, and pivot elements are set to zero as a result of applying the operator $\nu_{ij}$.

The form of the Jacobi annihilators is described by Theorem~\ref{tm:hjan} below. To this end set
$$\tau(s,t):=\left\{
              \begin{array}{cc}
                (t-1)(t-2)/2+s, & \text{for} \ 1\leq s<t\leq n, \\
                \tau(t,s)+N, & \text{for} \ 1\leq t<s\leq n. \\
              \end{array}
            \right.
$$
The function $\tau(s,t)$ specifies the position of the element $a_{st}$ in the vector $\ve(A)$.

\begin{Theorem}\label{tm:hjan}
Let $\mathcal{R}=\mathcal{R}_{ij}(U)$, $U=R(i,j,\phi,\alpha)$, be a Jacobi annihilator.
Then $\mathcal{R}$ differs from the identity matrix $I_{2N}$ in exactly $2n-2$ submatrices determined by
\[
\mathcal{R}_{\tau(i,j),\tau(i,j)}=0, \qquad \mathcal{R}_{\tau(j,i),\tau(j,i)}=0
\]
and
\begin{align*}
\left[
    \begin{array}{cc}
      \mathcal{R}_{\tau(r,i),\tau(r,i)} & \mathcal{R}_{\tau(r,i),\tau(r,j)} \\
      \mathcal{R}_{\tau(r,j),\tau(r,i)} & \mathcal{R}_{\tau(r,j),\tau(r,j)} \\
    \end{array}
  \right] & =\left[
            \begin{array}{cc}
              \cos\phi & e^{- \imath\alpha}sin\phi \\
              -e^{ \imath\alpha}sin\phi & \cos\phi \\
                         \end{array}
          \right], \\
\left[
    \begin{array}{cc}
      \mathcal{R}_{\tau(i,r),\tau(i,r)} & \mathcal{R}_{\tau(i,r),\tau(j,r)} \\
      \mathcal{R}_{\tau(j,r),\tau(i,r)} & \mathcal{R}_{\tau(j,r),\tau(j,r)} \\
    \end{array}
  \right] & =\left[
            \begin{array}{cc}
              \cos\phi & e^{ \imath\alpha}sin\phi \\
              -e^{- \imath\alpha}sin\phi & \cos\phi \\
                         \end{array}
          \right],
\end{align*}
wherer $1\leq r\leq n$, $r\notin\{i,j\}$.
\end{Theorem}

\begin{proof}
Since $\mathcal{R}$ annihilates the elements on positions $(i,j)$ and $(j,i)$, corresponding diagonal entries of $\mathcal{R}$ equal zero, i.e.\,
$$\mathcal{R}_{\tau(i,j),\tau(i,j)}=0, \quad \text{i} \quad \mathcal{R}_{\tau(j,i),\tau(j,i)}=0.$$
To get the other entries of $\mathcal{R}$ we use the definition of the Jacobi annihilator. Let of $A'=U^*AU$. For $i<j$ and $r\notin\{i,j\}$ have
\begin{align}
a_{ri}' & =a_{ri}\cos\phi+a_{rj}e^{- \imath\alpha}\sin\phi, &
a_{rj}' & =-a_{ri}e^{ \imath\alpha}\sin\phi+a_{rj}\cos\phi, \label{hjan1} \\
a_{ir}' & =a_{ir}\cos\phi+a_{jr}e^{ \imath\alpha}\sin\phi, &
a_{jr}' & =-a_{ir}e^{- \imath\alpha}\sin\phi+a_{jr}\cos\phi. \label{hjan2}
\end{align}
From relation~\eqref{hjan1} we obtain
$$\left[
    \begin{array}{c}
      a_{ri}' \\
      a_{rj}' \\
    \end{array}
  \right]=\left[
            \begin{array}{cc}
              \cos\phi & e^{- \imath\alpha}\sin\phi \\
              -e^{ \imath\alpha}\sin\phi & \cos\phi \\
            \end{array}
          \right]\left[
    \begin{array}{c}
      a_{ri} \\
      a_{rj} \\
    \end{array}
  \right].
$$
Hence,
$$\left[
    \begin{array}{cc}
      \mathcal{R}_{\tau(r,i),\tau(r,i)} & \mathcal{R}_{\tau(r,i),\tau(r,j)} \\
      \mathcal{R}_{\tau(r,j),\tau(r,i)} & \mathcal{R}_{\tau(r,j),\tau(r,j)} \\
    \end{array}
  \right]=\left[
            \begin{array}{cc}
              \cos\phi & e^{- \imath\alpha}\sin\phi \\
              -e^{ \imath\alpha}\sin\phi & \cos\phi \\
            \end{array}
          \right].$$
From~\eqref{hjan2} we have
$$\left[
    \begin{array}{c}
      a_{ir}' \\
      a_{jr}' \\
    \end{array}
  \right]=\left[
            \begin{array}{cc}
              \cos\phi & e^{ \imath\alpha}\sin\phi \\
              -e^{- \imath\alpha}\sin\phi & \cos\phi \\
            \end{array}
          \right]\left[
    \begin{array}{c}
      a_{ir} \\
      a_{jr} \\
    \end{array}
  \right],
$$
and consequently
$$\left[
    \begin{array}{cc}
      \mathcal{R}_{\tau(i,r),\tau(i,r)} & \mathcal{R}_{\tau(i,r),\tau(j,r)} \\
      \mathcal{R}_{\tau(j,r),\tau(i,r)} & \mathcal{R}_{\tau(j,r),\tau(j,r)} \\
    \end{array}
  \right]=\left[
            \begin{array}{cc}
              \cos\phi & e^{ \imath\alpha}\sin\phi \\
              -e^{- \imath\alpha}\sin\phi & \cos\phi \\
            \end{array}
          \right].$$
\end{proof}

In the sequel let $A$ be a Hermitian matrix. Let us consider how this affects the definition of Jacobi annihilators.
For Hermitian $A$, vector $\ve (A)$ has the form~\eqref{herm_A}, which means that $a=\ve(A)$ is an element of the real vector space
\[
\mathcal{H}_n = \left\{\left[\begin{array}{c}v \\ \bar{v} \end{array}\right]\ \Big{|} \ v\in \mathbb{C}^{N}\right\}.
\]
Thus, we can look at $\ve $ as a linear operator which maps the real vector space of Hermitian matrices of order $n$ onto the real vector space $\mathcal{H}_n$. Then $a\mapsto \mathcal{R}_{ij}(U)a$ will be the restriction to $\mathcal{H}_n$ of the general transformation from Definition~\ref{def: JacAnn}. Note that the matrix $\mathcal{R}_{ij}(U)$ will not be affected by the fact that $a\in \mathcal{H}_n$. Obviously, if $a\in \mathcal{H}_n$ then $\ve_0^{-1}(a)$ will be Hermitian matrix with zero diagonal elements.

Jacobi annihilators are building blocks of the Jacobi operators. Application of a Jacobi annihilator to vector $a=\ve(A)$ corresponds to a step of the Jacobi method on matrix $A$ and application of a Jacobi operator to $a=\ve(A)$ corresponds to a sweep of the cyclic Jacobi method applied to $A$.

\begin{Definition}
Let $N=n(n-1)/2$, $\nu\in [0,1]$, and $\mathcal{O} = (i_0,j_0),(i_1,j_1),\ldots,(i_{N-1},j_{N-1})$ $\in$ $\Ol(\mathcal{P}_n)$. Then
\[
\Jl_{\!\!\! \mathcal{O}}^{\nu} = \{\mathcal{J} \ \big{|} \ \mathcal{J}=\mathcal{R}_{i_{N-1}j_{N-1}}\ldots\mathcal{R}_{i_1j_1}\mathcal{R}_{i_0j_0},
\ \mathcal{R}_{i_kj_k}\in\Rl_{i_k j_k}^{\nu}, \ 0\leq k\leq N-1\}
\]
is the class of Jacobi operators associated with the sequence $\mathcal{O}$ and $\nu$. The $2N\times2N$ matrices $\mathcal{J}$ from the $\Jl_{\!\!\! \mathcal{O}}^{\nu}$ are called the Jacobi operators.
If $\nu =0$, we write $\Jl_{\!\!\! \mathcal{O}}$ instead of $\Jl_{\!\!\! \mathcal{O}}^{0}$.
\end{Definition}

Recall that $\Ol(\mathcal{P}_n)$ is the set of all orderings of $\mathcal{P}_n$ and these orderings define all cyclic strategies. If we know the relation between two orderings, we can say something about the corresponding Jacobi operators.

\begin{Proposition}\label{prop:2.6}
Let $\mathcal{O}\in\Ol(\mathcal{P}_n)$, $\nu\in [0,1]$, and suppose $\|\mathcal{J}\|_2\leq\mu$ for all $\mathcal{J}\in\Jl_{\!\!\! \mathcal{O}}^{\nu}$.
\begin{itemize}
\item[(i)]\quad\ If $\mathcal{O}'=\mathcal{O}(\mathsf{p})$ and $\mathcal{J}'\in\Jl_{\!\!\! \mathcal{O}'}^{\nu}$, then $\|\mathcal{J}'\|_2\leq\mu$.
\item[(ii)]\quad\  If $\mathcal{J}'' \in \Jl_{\!\!\! \mathcal{O}^{\leftarrow}}^{\nu}$, then $\|\mathcal{J}''\|_2\leq\mu$.
\end{itemize}
\end{Proposition}

\begin{proof}
The proofs are identical to those from~\cite{har+beg-17} which deal with real matrices and real Jacobi annihilators. 
Assertion $(i)$ is a special case of~\cite[Theorem~2.22(i)]{har+beg-17} and assertion $(ii)$ is a special case of~\cite[Proposition~2.21]{har+beg-17}. We have to use the special partition $\pi =(1,1,\ldots ,1)$.
The parameter $\nu$ corresponds to $\rho/\sqrt{2}$ in
\cite{har+beg-17}, and the proofs hold for any $0\leq \nu\leq 1$.
\end{proof}

\begin{Proposition}\label{tm:joppinv}
Let $\mathcal{O}\in\Ol(\mathcal{P}_n)$, $\nu\in [0,1]$, and suppose $\|\mathcal{J}\|_2\leq\mu$ for all $\mathcal{J}\in\Jl_{\!\!\! \mathcal{O}}^{\nu}$.
If $\mathcal{O}'\stackrel{\mathsf{w}}{\sim}\mathcal{O}$ and there are exactly $d$ relations of shift equivalence in the chain connecting $\mathcal{O}$ and $\mathcal{O}'$,
then for any $d+1$ Jacobi operators from $\Jl_{\!\!\! \mathcal{O}'}^{\nu}$ we have
\[
\|\mathcal{J}_1' \mathcal{J}_2'\cdots \mathcal{J}_{d+1}'\|_2 \leq \mu, \quad \mathcal{J}_1',\ldots,\mathcal{J}_{d+1}' \in \Jl_{\!\!\! \mathcal{O}'}^{\nu}.
\]
\end{Proposition}

\begin{proof} This is a special case of~\cite[Proposition~2.20]{har-15}. \end{proof}

\section{Convergence of the complex Jacobi method for different pivot strategies}\label{sec:strategies}

The most common cyclic pivot strategies are the serial pivot strategies: row-cyclic strategy $I_{\text{row}}=I_{\mathcal{O}_{\text{row}}}$ and column-cyclic strategy $I_{\text{col}}=I_{\mathcal{O}_{\text{col}}}$. They are defined by the ``row-wise'' and ``column-wise'' orderings of $\mathcal{P}_n$,
\begin{align*}
\mathcal{O}_{\text{row}} & = (1,2),(1,3),\ldots,(1,n), (2,3),\ldots,(2,n),\ldots,(n-1,n), \\
\mathcal{O}_{\text{col}} & = (1,2),(1,3),(2,3),\ldots,(1,n), (2,n),\ldots,\ldots,(n-1,n).
\end{align*}
Cyclic pivot strategies equivalent (resp.\@ weak-equivalent) to the serial ones are called wavefront (resp.\@ weak-wavefront) strategies (see~\cite{SS89}).
The class of  serial pivot strategies with permutations introduced in~\cite{B,har+beg-17} is derived from the class of serial strategies.

Denote the set of all permutations of the set $\{l_1,l_1+1,\ldots,l_2\}$ by $\Pi^{(l_1,l_2)}$. Let
\begin{align*}
\mathcal{C}_c^{(n)} = \big{\{} \mathcal{O}\in\mathcal{\Ol}(\mathcal{P}_n) \ \big{|} \ \mathcal{O}= & (1,2),(\tau_{3}(1),3),(\tau_{3}(2),3),\ldots,(\tau_{n}(1),n),\ldots \\
& \ldots,(\tau_{n}(n-1),n), \quad \tau_{j}\in\Pi^{(1,j-1)}, \ 3\leq j\leq n \big{\}}. \nonumber
\end{align*}
The orderings from $\mathcal{C}_c^{(n)}$ are called \emph{column-wise orderings with permutations} of the set $\mathcal{P}_n$.
First pair of $\mathcal{O}\in\mathcal{C}_c^{(n)}$ is $(1,2)$, then we have all pairs from the second column in some order, then the pairs from the third column, etc. At the last stage we have all pairs from the $n$th column in some order.

\begin{Theorem}\label{tm:c1}
Let $A$ be a Hermitian matrix of dimension $n$ and $\mathcal{O}\in\mathcal{C}_c^{(n)}$. Suppose that $A'$ is obtained from $A$ by applying one cycle of the cyclic Jacobi method defined by the strategy $I_{\mathcal{O}}$.
If all rotation angles satisfy $|\cos\phi_k| \geq \nu>0$, $k\geq0$,
then there is a constant $\gamma_{n,\nu}$ depending only on $n$ and $\nu$ such that
\begin{equation}\label{rel:c1}
S^2(A')\leq \gamma_{n,\nu} S^2(A), \quad 0\leq \gamma_{n,\nu}<1.
\end{equation}
\end{Theorem}

\begin{proof}
We use mathematical induction over the matrix size $n$. For $n=2$, it is easy to see that relation~\eqref{rel:c1} holds with $\gamma_{2,\nu}=0$.
Assume that~\eqref{rel:c1} holds for $n-1$, with $0\leq \gamma_{n-1,\nu}<1$.

Denote by $\tilde{A}=(\tilde{a}_{ij})$ a matrix obtained from $A$ after the first $\tilde{N}=(n-1)(n-2)/2$ steps of the Jacobi method.
Let us observe the last column of $\tilde{A}$. Each of $\tilde{N}$ Jacobi rotations used to get from $A$ to $\tilde{A}$ changed only two elements from the last column. Thereby the sum of the squares of their absolute values remained the same,
\begin{equation}\label{c1suma}
\sum_{i=1}^{n-1}|a_{in}|^2=\sum_{i=1}^{n-1}|\tilde{a}_{in}|^2.
\end{equation}
We define $0\leq\epsilon\leq1$ such that
\begin{equation}\label{c1eps}
(1-\epsilon^2)S^2(A)=\sum_{i=1}^{n-1}|a_{in}|^2.
\end{equation}
By $A_{n-1}$ (resp.\@ $\tilde{A}_{n-1}$) denote the leading $(n-1)\times(n-1)$ submatrix of $A$ (resp.\@ $\tilde{A}$). From the induction hypothesis we have
\begin{equation}\label{c1Sn-1}
S^2(\tilde{A}_{n-1})\leq\gamma_{n-1,\nu}S^2(A_{n-1}).
\end{equation}
Using relations~\eqref{c1Sn-1}, \eqref{c1suma}  and~\eqref{c1eps} we obtain
\begin{align}
S^2(\tilde{A}) =S^2(\tilde{A}_{n-1})+2\sum_{i=1}^{n-1}|\tilde{a}_{in}|^2 &\leq \gamma_{n-1,\nu}S^2(A_{n-1})+2\sum_{i=1}^{n-1}|a_{in}|^2 \nonumber \\
& =\gamma_{n-1,\nu}(S^2(A)-(1-\epsilon^2)S^2(A))+(1-\epsilon^2)S^2(A) \nonumber \\
& \leq (1-\epsilon^2(1-\gamma_{n-1,\nu}))S^2(A). \label{c1Stilda}
\end{align}
We cannot set $\gamma_{n,\nu}=1-\epsilon^2(1-\gamma_{n-1,\nu})$ because $\epsilon$ can be arbitrary small and it could produce $\gamma_{n,\nu}$ equal to one. Thus, we need to analyse how the last $n-1$ rotations affect the off-norm of the underlaying matrix.

By $\tilde{A}^{(k)}=(\tilde{a}_{st}^{(k)})$ denote the matrix obtained from $\tilde{A}$ after $k$ additional rotations. This corresponds to $\tilde{k}=\tilde{N}+k$ rotations on $A$. For $1\leq k\leq n-1$, pivot element of $\tilde{A}^{(k)}$ is $\tilde{a}_{\pi_{n}(k),n}^{(k-1)}=a_{\pi_{n}(k),n}^{(\tilde{N}+k-1)}$.
Rotation angles $\phi_{\tilde{N}+k}$ and $\alpha_{\tilde{N}+k}$ are determined by the pivot element $\tilde{a}_{\pi_{n}(k),n}^{(k-1)}$.
Denote $c_k=\cos\phi_{\tilde{N}+k}$, $s_k=\sin\phi_{\tilde{N}+k}$, $1\leq k\leq n-1$.

We have to estimate $|\tilde{a}_{\pi_{n}(i),n}^{(i-1)}|^2$ for $1\leq i\leq n-1$. To obtain the expression for $\tilde{a}_{\pi_{n}(i),n}^{(i-1)}$ we use relation
\[
\tilde{a}_{\pi_{n}(i),n}^{(k)} = c_k\tilde{a}_{\pi_{n}(i),n}^{(k-1)} - e^{\imath\tilde{\alpha}_k}s_{k}\tilde{a}_{\pi_{n}(i),\pi_{n}(k)}^{(k-1)}, \quad 1\leq k\leq i-1.
\]
Elements in the $\pi_{n}(i)$th row of $\tilde{A}_{n-1}$ are changed exactly once during the last $n-1$ steps. Hence, $\tilde{a}_{\pi_{n}(i),\pi_{n}(k)}^{(k-1)}=\tilde{a}_{\pi_{n}(i),\pi_{n}(k)}$ for $1\leq k\leq n-1$. This implies
\begin{equation}\label{c1a}
\tilde{a}_{\pi_{n}(i),n}^{(k)} = c_{k}\tilde{a}_{\pi_{n}(i),n}^{(k-1)} - e^{\imath\tilde{\alpha}_k}s_{k}\tilde{a}_{\pi_{n}(i),\pi_{n}(k)}, \quad 1\leq k\leq i-1.
\end{equation}
Pivot element in the step $\tilde{N}+i$ is $\tilde{a}_{\pi_{n}(i),n}^{(i-1)}$.
We write down relation~\eqref{c1a} for the consecutive values of $k$ and get
\begin{align*}
\tilde{a}_{\pi_{n}(i),n}^{(1)} & = c_1\tilde{a}_{\pi_{n}(i),n}-s_1e^{\imath\tilde{\alpha}_1}\tilde{a}_{\pi_{n}(i),\pi_{n}(1)}\\
\tilde{a}_{\pi_{n}(i),n}^{(2)} & = c_2\tilde{a}_{\pi_{n}(i),n}^{(1)} -s_2e^{\imath\tilde{\alpha}_2}\tilde{a}_{\pi_{n}(i),\pi_{n}(2)}\\
&  \vdots \\
\tilde{a}_{\pi_{n}(i),n}^{(i-2)} & = c_{i-2}\tilde{a}_{\pi_{n}(i),n}^{(i-3)} -s_{i-2}e^{\imath\tilde{\alpha}_{i-2}}\tilde{a}_{\pi_{n}(i),\pi_{n}(i-2)} \\
\tilde{a}_{\pi_{n}(i),n}^{(i-1)} & = c_{i-1}\tilde{a}_{\pi_{n}(i),n}^{(i-2)} -s_{i-1}e^{\imath\tilde{\alpha}_{i-1}}\tilde{a}_{\pi_{n}(i),\pi_{n}(i-1)} .
\end{align*}
Note that $\tilde{a}_{\pi_{n}(i),n}^{(i)} = 0$.
To get $\tilde{a}_{\pi_{n}(i),n}^{(i-1)}$ we multiply the equation for $\tilde{a}_{\pi_{n}(i),n}^{(i-2)}$ by $c_{i-1}$, the equation for $\tilde{a}_{\pi_{n}(i),n}^{(i-3)}$ by $c_{i-1}c_{i-2}$, etc., until the equation for $\tilde{a}_{\pi_{n}(i),n}^{(1)}$ which we multiply by $c_{i-1}c_{i-2}\cdots c_2$. After the summation of newly obtained equations we get
\begin{equation}\label{c1a(i-1)}
\tilde{a}_{\pi_{n}(i),n}^{(i-1)} = c_1\ldots c_{i-1}\tilde{a}_{\pi_{n}(i),n} - \sum_{k=1}^{i-1} s_kc_{k+1}\ldots c_{i-1}e^{ i \tilde{\alpha}_k}\tilde{a}_{\pi_{n}(i),\pi_{n}(k)},
\end{equation}
where for $k=i-1$ the product $s_kc_{k+1}\ldots c_{i-1}$ reduces to $s_{i-1}$. Relation~\eqref{c1a(i-1)} holds for $1\leq i\leq n-1$ and for $i=1$ it reduces to $\tilde{a}_{\pi_{n}(i),n}=\tilde{a}_{\pi_{n}(i),n}$.

Until the end of the cycle reduction of $S^2(\tilde{A})$ equals $2\sum_{i=1}^{n-1}\big{|}\tilde{a}_{\pi_{n}(i),n}^{(i-1)}\big{|}^2$. We need to determine its lower bound.
To bound the value $\big{|}\tilde{a}_{\pi_{n}(i),n}^{(i-1)}\big{|}^2$ we use the Cauchy-Schwarz inequality and the inequality $|a-b|^2\geq\frac{1}{2}|a|^2-|b|^2$ that holds for all $a,b\in\mathbb{C}$.
Since $|\cos\phi_i|\geq\nu$, we have $c_i^2\geq\nu^2$ and $s_i^2\leq1-\nu^2$. Moreover, we have $|e^{\pm \imath \tilde{\alpha}_k}|=1$ for all $k$. Thus, from~\eqref{c1a(i-1)} it follows
\begin{align*}
\big{|}\tilde{a}_{\pi_{n}(i),n}^{(i-1)}\big{|}^2 & \geq \frac{1}{2}\big{|}c_1\cdots c_{i-1}\tilde{a}_{\pi_{n}(i),n}\big{|}^2
- \Big{|}\sum_{k=1}^{i-1} s_kc_{k+1}\cdots c_{i-1}\tilde{a}_{\pi_{n}(i),\pi_{n}(k)}\Big{|}^2\\
& \geq \frac{1}{2}c_1^2c_2^2\cdots c_{i-1}^2 \big{|}\tilde{a}_{\pi_{n}(i),n}\big{|}^2 -
\Big{(}\sum_{k=1}^{i-1}s_k^2c_{k+1}^2\cdots c_{i-1}^2\Big{)}
\Big{(}\sum_{k=1}^{i-1}|\tilde{a}_{\pi_{n}(i),\pi_{n}(k)}|^2\Big{)},
\end{align*}
for $1\leq i\leq n-1$. It is easy to show that
\[
s_{i-1}^2+s_{i-2}^2c_{i-1}^2+s_{i-3}^2c_{i-2}^2c_{i-1}^2+\cdots\ +s_{1}^2c_{2}^2\cdots c_{i-1}^2 =
1- c_1^2c_2^2\cdots c_{i-1}^2.
\]
Therefore, we have
\[
\big{|}\tilde{a}_{\pi_{n}(i),n}^{(i-1)}\big{|}^2  \geq
\frac{1}{2}c_1^2c_2^2\cdots c_{i-1}^2 \big{|}\tilde{a}_{\pi_{n}(i),n}\big{|}^2 -
(1- c_1^2c_2^2\cdots c_{i-1}^2)\sum_{k=1}^{i-1}|\tilde{a}_{\pi_{n}(i),\pi_{n}(k)}|^2.
\]
We use the following inequalities
\begin{align*}
& c_1^2c_2^2\cdots c_{i-1}^2\geq \nu^{2(i-1)}\geq \nu^{2(n-2)}, \\
& 1-c_1^2c_2^2\cdots c_{i-1}^2 \leq 1-\nu^{2(i-1)}\leq 1-\nu^{2(n-2)},
\end{align*}
that hold for $1\leq i\leq n-1$.
Now, making the sum over $1\leq i\leq n-1$ yields
\begin{align*}
2\sum_{i=1}^{n-1}\big{|}\tilde{a}_{\pi_{n}(i),n}^{(i-1)}\big{|}^2 & \geq
\sum_{i=1}^{n-1} \nu^{2(i-1)}\big{|}\tilde{a}_{\pi_{n}(i),n}\big{|}^2 -
2\sum_{i=1}^{n-1}(1-\nu^{2(i-1)})\sum_{k=1}^{i-1}|\tilde{a}_{\pi_{n}(i),\pi_{n}(k)}|^2 \\
 & \geq \nu^{2(n-2)}\sum_{i=1}^{n-1} \big{|}\tilde{a}_{\pi_{n}(i),n}\big{|}^2 - (1-\nu^{2(n-2)})\,2\sum_{i=1}^{n-1}\sum_{k=1}^{i-1}|\tilde{a}_{\pi_{n}(i),\pi_{n}(k)}|^2 .
\end{align*}
To bound the first term on the right-hand side we use~\eqref{c1eps}. We have
\[
\sum_{i=1}^{n-1}\big{|}\tilde{a}_{\pi_{n}(i),n}\big{|}^2  =\sum_{i=1}^{n-1}\big{|}\tilde{a}_{in}\big{|}^2
=\sum_{i=1}^{n-1}\big{|}a_{in}\big{|}^2 = (1-\epsilon^2)\frac{1}{2}S^2(A).
\]
To bound the second term, we use permutation matrix $P$ of dimension $n-1$, defined by the permutation $\pi_n$ of the set $\{1,2,\ldots,n-1\}$.
For $I_{n-1}=[e_1,\ldots,e_{n-1}]$ we have
\[
P e_j =e_{\pi_{n}(j)}, \quad 1\leq j\leq n-1.
\]
Since
\[
(P^T\tilde{A}_{n-1}P)_{ik} = e_i^T P^T \tilde{A}_{n-1} Pe_k =
e_{\pi_{n}(i)}^T \tilde{A}_{n-1}e_{\pi_{n}(k)} =\tilde{a}_{\pi_{n}(i),\pi_{n}(k)}
\]
for $ 1\leq i,k\leq n-1$, we conclude that the transformation $\tilde{A}_{n-1}\mapsto P^T\tilde{A}_{n-1}P$ leaves the diagonal elements on the diagonal. Therefore, we have $S^2(P^T\tilde{A}_{n-1}P)=S^2(\tilde{A}_{n-1})$ and using relation~\eqref{c1Sn-1} we obtain
\begin{align*}
2\sum_{i=1}^{n-1}\sum_{k=1}^{i-1}|\tilde{a}_{\pi_{n}(i),\pi_{n}(k)}|^2 &=
2\sum_{i=1}^{n-1}\sum_{k=1}^{i-1}|(P^T\tilde{A}_{n-1}P)_{ik}|^2 = S^2(P^T\tilde{A}_{n-1}P) =
S^2(\tilde{A}_{n-1})\\
&\leq  \gamma_{n-1,\nu}S^2(A_{n-1}) =  \gamma_{n-1,\nu}\epsilon^2 S^2(A).
\end{align*}

From the above estimates it follows that
\begin{equation}\label{c1nejedn}
2\sum_{i=1}^{n-1}\big{|}\tilde{a}_{\pi_{n}(i),n}^{(i-1)}\big{|}^2
\geq f_{n,\nu}(\epsilon)S^2(A),
\end{equation}
where the function $f_{n,\nu}:[0,1]\rightarrow\mathbb{R}$ is defined by
\[
f_{n,\nu}(\epsilon) = \frac{1}{2}\nu^{2(n-2)}-\Big{(}\frac{1}{2}\nu^{2(n-2)}+ \gamma_{n-1,\nu}\big{(}1-\nu^{2(n-2)}\big{)}\Big{)}\epsilon^2.
\]
This function is decreasing, attaining its maximum $\frac{\nu^{2(n-2)}}{2}$ at $0$, minimum
$-\gamma_{n-1,\nu}(1-\nu^{2(n-2)})$ at $1$, and zero at
\[
\epsilon_{n,\nu}^{(0)} = \frac{\nu^{n-2}}{\sqrt{\nu^{2(n-2)}+2\gamma_{n-1,\nu}(1-\nu^{2(n-2)})}}.
\]
The left-hand side in~\eqref{c1nejedn} is nonnegative, thus $f_{n,\nu}$ can be replaced by the nonnegative function
\begin{equation}\label{c1f+}
f_{n,\nu}^{+}(\epsilon)=\left\{
                    \begin{array}{ll}
                      f_{n,\nu}(\epsilon), & \epsilon\in[0,\epsilon_{n,\nu}^{(0)}], \\
                      0, & \epsilon\in[\epsilon_{n,\nu}^{(0)},1].
                    \end{array}
                  \right.
\end{equation}
Function $f_{n,\nu}^{+}$ is continuous. Therefore, there is a real number  $\epsilon_{n,\nu}\in \langle 0,\epsilon_{n,\nu}^{(0)}\rangle$ such that
\[
f_{n,\nu}^{+}(\epsilon_{n,\nu}) = \frac{1}{4}\nu^{2(n-2)}.
\]

Now we can bound $S^2(A')$. Using the relations~\eqref{c1Stilda}, \eqref{c1nejedn}, \eqref{c1f+} and~\eqref{c1Stilda} we get
\[
S^2(A') = S^2(\tilde{A})-2\sum_{i=1}^{n-1}\big{|}\tilde{a}_{\pi_{n-1}(i),n}^{(i-1)}\big{|}^2 \leq g_{n,\nu}(\epsilon)S^2(A),
\]
where the function $g_{n,\nu}:[0,1]\rightarrow\mathbb{R}$ is given by
\[
g_{n,\nu}(\epsilon)=1-\epsilon^2(1-\gamma_{n-1,\nu})-f_{n,\nu}^{+}(\epsilon),\quad \epsilon\in[0,1].
\]
We obviously have
\[
g_{n,\nu}(\epsilon) \leq\left\{
                    \begin{array}{ll}
                      1-\frac{1}{4}\nu^{2(n-2)}, & \text{if }\epsilon\in[0,\epsilon_{n,\nu}], \\
                      1-\epsilon^2(1-\gamma_{n-1,\nu}), & \text{if }\epsilon\in[\epsilon_{n,\nu},1].
                    \end{array}
                  \right.
\]
Hence, we can end the proof by setting
\[
\gamma_{n,\nu} = \max \left\{1-\frac{1}{4}\nu^{2(n-2)},1-\epsilon_{n,\nu}^2(1-\gamma_{n-1,\nu})\right\} .
\]
\end{proof}

The result from the previous theorem can be obtained for the more general process with the Jacobi operators.

\begin{Theorem}\label{tm:c1jop}
Let $\mathcal{O}\in\mathcal{C}_c^{(n)}$, $\nu\in (0,1]$, and $\mathcal{J}\in\Jl_{\!\!\! \mathcal{O}}^{\nu}$. Then there is a constant $\mu_{n,\nu}$ depending only on $n$ and $\nu$ such that
\[
\|\mathcal{J}\|_2\leq\mu_{n,\nu}, \quad 0\leq\mu_{n,\nu}<1.
\]
\end{Theorem}

\begin{proof}
Let $a\in\mathcal{H}_n$ be an arbitrary vector and let $a'=\mathcal{J}a$, that is
\begin{equation}\label{c1jopa}
a'=\mathcal{J}a = \mathcal{R}_{\pi_n(n-1),n}\cdots\mathcal{R}_{\pi_3(2),3}\mathcal{R}_{\pi_3(1),3}\mathcal{R}_{1,2} \ a, \quad \mathcal{R}_{ij}\in\Rl_{ij}^{\nu}.
\end{equation}
Matrices $A=\ve_0^{-1}(a), A'=\ve_0^{-1}(a')\in S_0$ are two-dimensional representations of the vectors $a$ and $a'$, respectively. They both are Hermitian and we have
\[
S^2(A)=\|a\|_2^2, \qquad S^2(A')=\|a'\|_2^2=\|\mathcal{J} a\|_2^2.
\]
The process~\eqref{c1jopa} can be written in the form
\[
a^{(0)}=a,\ a^{(1)}= \mathcal{R}_{1,2}a^{(0)},\ a^{(2)}= \mathcal{R}_{\pi_3(1),3}a^{(1)},\ \ldots\ ,
a^{(n-1)}= \mathcal{R}_{\pi_n(n-1),n}a^{(n-2)}
\]
and we have $a'=a^{(n-1)}$. We can display vectors $a^{(k)}$, $0\leq k\leq n-1$, as two-dimensional arrays, which are Hermitian matrices with zero diagonal elements. Comparing aforesaid sequence of Hermitian matrices to the standard complex Jacobi method ~\eqref{hjacobiagm} we can notice that both amount to almost the same procedure.
Although they generate different iteration matrices, since the rotation angles are different, all estimates used in the proof of Theorem~\ref{tm:c1} still hold. That is why we can conclude that
\[
S^2(A')\leq \gamma_{n,\nu} S^2(A),
\]
where $0\leq\gamma_{n,\nu}<1$ is the same as in Theorem~\ref{tm:c1}.
Therefore we have $\|\mathcal{J}a\|_2^2\leq\gamma_{n,\nu} \|a\|_2^2$, and
\[
\|\mathcal{J}\|_2=\max_{a\in \mathcal{H}_n\setminus \{0\} }\frac{\|\mathcal{J}a\|_2}{\|a\|_2}\leq\mu_{n,\nu}, \quad \mu_{n,\nu}=\sqrt{\gamma_{n,\nu}}.
\]
\end{proof}

\subsection{Serial strategies with permutations}

Similarly to the set of column-wise ordering with permutations, we define the set of \emph{row-wise orderings with permutations}.
Let
\begin{align*}
\mathcal{C}_r^{(n)} = \big{\{}\mathcal{O}\in\mathcal{\Ol}(\mathcal{P}_n) \ \big{|} & \ \mathcal{O}=(n-1,n),(n-2,\tau_{n-2}(n-1)),(n-2,\tau_{n-2}(n)),\ldots \\
& \qquad \ldots,(1,\tau_{1}(2)),\ldots,(1,\tau_{1}(n)), \quad \tau_{i}\in\Pi^{(i+1,n)}, \ 1\leq i\leq n-2 \big{\}}.
\end{align*}
Here we first have pair $(n-1,n)$, then in some order all pairs from the $(n-2)$nd row, then all pairs from the $(n-3)$rd row, etc. In the last stage we have the pairs from the first row in some order. It is easy to check that every ordering from $\mathcal{C}_r^{(n)}$ is permutation equivalent to an ordering from $\mathcal{C}_c^{(n)}$ with permutation
$\left(
   \begin{array}{cccc}
     1 & 2 & \ldots & n \\
     n & n-1 & \ldots & 1 \\
   \end{array}
 \right)
$.
From the sets of orderings $\mathcal{C}_c^{(n)}$ and $\mathcal{C}_r^{(n)}$ we can easily construct other two sets, equally large and closely related. These are the sets of orderings inverse to those from $\mathcal{C}_c^{(n)}$ and $\mathcal{C}_r^{(n)}$,
$$\overleftarrow{\mathcal{C}}_c^{(n)} = \big{\{} \mathcal{O}\in\mathcal{\Ol}(\mathcal{P}_n) \ \big{|} \ \mathcal{O}^{\leftarrow }\in \mathcal{C}_c^{(n)} \big{\}} \qquad \text{and} \qquad
\overleftarrow{\mathcal{C}}_r^{(n)} = \big{\{} \mathcal{O}\in\mathcal{\Ol}(\mathcal{P}_n) \ \big{|} \ \mathcal{O}^{\leftarrow }\in \mathcal{C}_r^{(n)} \big{\}}.$$
These four sets of pivot orderings form the set of \emph{serial orderings with permutations},
\begin{equation}\label{def:sp}
\mathcal{C}_{sp}^{(n)} = \mathcal{C}_c^{(n)}\cup\overleftarrow{\mathcal{C}}_c^{(n)} \cup \mathcal{C}_r^{(n)}\cup\overleftarrow{\mathcal{C}}_r^{(n)}.
\end{equation}

Theorems~\ref{tm:c1} and~\ref{tm:c1jop} hold if we replace $\mathcal{C}_c^{(n)}$ with $\mathcal{C}_{sp}^{(n)}$.

\begin{Theorem}\label{tm:spjop}
Let $\mathcal{O}\in\mathcal{C}_{sp}^{(n)}$, $\nu\in (0,1]$ and let $\mathcal{J}\in\Jl_{\!\!\! \mathcal{O}}^{\nu}$. Then there is constant $\mu_{n,\nu}$ depending only on $n$ and $\nu$ such that
\[
\|\mathcal{J}\|_2\leq\mu_{n,\nu}, \quad 0\leq\mu_{n,\nu}<1.
\]
\end{Theorem}

\begin{proof}
If $\mathcal{O}\in\mathcal{C}_c^{(n)}$, this boils down to Theorem~\ref{tm:c1jop}.
If $\mathcal{O}\in\mathcal{C}_r^{(n)}$, we use Theorem~\ref{tm:c1jop} and Proposition~\ref{prop:2.6}(i).
If $\mathcal{O}\in\overleftarrow{\mathcal{C}}_c^{(n)} \cup \overleftarrow{\mathcal{C}}_r^{(n)}$, we use Theorem~\ref{tm:c1jop} and Proposition~\ref{prop:2.6}(ii).
\end{proof}

\begin{Corollary}\label{tm:sp}
Let $A$ be a Hermitian matrix of order $n$ and $\mathcal{O}\in\mathcal{C}_{sp}^{(n)}$. Let $A'$ be obtained from $A$ by applying one sweep of the cyclic Jacobi method defined by $I_{\mathcal{O}}$.
If all rotation angles satisfy $\phi(k)\in[-\frac{\pi}{4},\frac{\pi}{4}]$, $k\geq0$,
then there is a constant $\gamma_n$ depending only on $n$ such that
$$S^2(A')\leq \gamma_n S^2(A), \quad 0\leq \gamma_n<1.$$
\end{Corollary}

\begin{proof}
Set $a=\ve(A)$ and $a'=\ve(A')$. Then $a,\,a'\in \mathcal{H}_n$ and we have $a'=\mathcal{J}_{\mathcal{O}} \, a$ for some $\mathcal{J}_{\mathcal{O}}$ from $\Jl_{\!\!\! \mathcal{O}}^{\pi/4}$. The Jacobi annihilators in $\mathcal{J}_{\mathcal{O}}$ are defined by the complex Jacobi rotations appearing in the considered sweep of the Jacobi method. From Theorem~\ref{tm:spjop} we obtain $\|\mathcal{J}_{\mathcal{O}}\|_2 \leq \mu_{n,\pi/4}$,  $0\leq \mu_{n,\pi/4}<1$. Thus we have
\[
S^2(A')=\|a'\|_2^2 = \|\mathcal{J}_{\mathcal{O}}a\|_2^2\leq \mu_{n,\pi/4}^2 \|a\|^2 = \gamma_n S^2(A),\quad \gamma_n =\mu_{n,\pi/4}^2.
\]
\end{proof}

\subsection{Generalized serial pivot strategies with permutations}

Orderings from~\eqref{def:sp} can be further generalized using the relations of weak and permutational equivalence (cf. \cite{har+beg-17}).

\begin{Definition}\label{def:sg}
The class of generalized serial orderings of $\mathcal{P}_n$ is given with
\[
\mathcal{C}_{sg}^{(n)} = \big{\{}\mathcal{O}\in\Ol(\mathcal{P}_n) \ \big{|} \ \mathcal{O}\stackrel{\mathsf{p}}{\sim}\mathcal{O'}\stackrel{\mathsf{w}}{\sim}\mathcal{O''} \ \text{or} \ \mathcal{O}\stackrel{\mathsf{w}}{\sim}\mathcal{O'}\stackrel{\mathsf{p}}{\sim}\mathcal{O''}, \ \mathcal{O''}\in\mathcal{C}_{sp}^{(n)}\big{\}},
\]
where weak equivalence relations are given in canonical form and $\mathcal{O'}\in\Ol(\mathcal{P}_n)$.
\end{Definition}

Compared to the Theorem~\ref{tm:spjop} and Corollary~\ref{tm:sp}, slightly weaker results hold if we use the orderings from $\mathcal{C}_{sg}^{(n)}$ instead.

\begin{Theorem}\label{tm:sg}
Let $\mathcal{O}\in\mathcal{C}_{sg}^{(n)}$. Suppose that $\mathcal{O}\stackrel{\mathsf{p}}{\sim}\mathcal{O'}\stackrel{\mathsf{w}}{\sim}\mathcal{O''}$  or
$\mathcal{O}\stackrel{\mathsf{w}}{\sim}\mathcal{O'}\stackrel{\mathsf{p}}{\sim}\mathcal{O''}$,
$\mathcal{O''}\in\mathcal{C}_{sp}^{(n)}$ and that the weak equivalence relation is in canonical form
containing $d$ shift equivalences.
Then for any $d+1$ Jacobi operators $\mathcal{J}_1,\mathcal{J}_2,\ldots,\mathcal{J}_{d+1}\in\Jl_{\!\!\! \mathcal{O}}^{\nu}$, $0<\nu\leq 1$, there is a constant $\zeta_{n,\nu}$ depending only on $n$ and $\nu$ such that for any $d+1$ Jacobi operators $\mathcal{J}_1,\mathcal{J}_2,\ldots,\mathcal{J}_{d+1}\in\Jl_{\!\!\! \mathcal{O}}^{\nu}$ one has
\[
\|\mathcal{J}_1\mathcal{J}_2\cdots\mathcal{J}_{d+1}\|_2 \leq\zeta_{n,\nu}, \quad 0\leq\zeta_{n,\nu}<1.
\]
\end{Theorem}

\begin{proof} This is a special case of~\cite[Theorem~3.9]{har+beg-17}. In that proof instead of \cite[Theorem~3.5]{har+beg-17}, \cite[Theorem~2.22(i)]{har+beg-17} and \cite[Proposition~2.20]{har+beg-17}, one can use Theorem~\ref{tm:spjop}, Proposition~\ref{prop:2.6}(i) and Proposition~\ref{tm:joppinv}, respectively.
\end{proof}

\begin{Corollary}\label{cor:3.7}
Let $A$ be a Hermitian matrix of order $n$.
Let $\mathcal{O}\in\mathcal{C}_{sg}^{(n)}$.
Suppose that $\mathcal{O}\stackrel{\mathsf{p}}{\sim}\mathcal{O'}\stackrel{\mathsf{w}}{\sim}\mathcal{O''}$  or
$\mathcal{O}\stackrel{\mathsf{w}}{\sim}\mathcal{O'}\stackrel{\mathsf{p}}{\sim}\mathcal{O''}$,
$\mathcal{O''}\in\mathcal{C}_{sp}^{(n)}$ and that the weak equivalence relation is in canonical form
containing $d$ shift equivalences.
Let $A'$ be obtained from $A$ by applying $d+1$ cycles of the cyclic Jacobi method defined by the strategy $I_{\mathcal{O}}$.
If all rotation angles satisfy $\phi(k)\in[-\frac{\pi}{4},\frac{\pi}{4}]$, $k\geq0$, then there is a constant $\gamma_n$ depending only on $n$ such that
\[
S^2(A')\leq \gamma_n S^2(A), \quad 0\leq \gamma_n<1.
\]
\end{Corollary}

\begin{proof} This is a special case of~\cite[Theorem~3.11]{har+beg-17} and also of Theorem~\ref{tm:general} from the next subsection.  The proof is similar to the proof of Corollary~\ref{tm:sp}. However it uses $d+1$ Jacobi operators instead of $N$ Jacobi annihilators.
\end{proof}
As in \cite{har+beg-17} it can be shown that if the chains in Definition~\ref{def:sg} do not use shifts, then for such orderings Theorem~\ref{tm:spjop} and Theorem~\ref{tm:sp} hold.

\subsection{A Jacobi-type process for Hermitian matrices}

In this subsection we formulate and prove the theorem on the convergence of element-wise Jacobi-type processes for Hermitian matrices. After we have proved the bounds for the complex Jacobi operators in previous subsections, the proof of this more general result becomes similar to the corresponding proofs in~\cite[Corollary~5.8]{har-15} and~\cite[Corollary~5.3]{har+beg-17}.

Elementary plane matrix is a nonsingular matrix that differs from the identity matrix $I_n$ in one principal submatrix of order $2$.
This submatrix is called pivot submatrix. In the rest of the paper we denote pivot indices by $i(k)$, $j(k)$, instead of $i_k$, $j_k$, respectively.

\begin{Theorem}\label{tm:general}
Let $H\neq0$ be a Hermitian matrix. Let $(H^{(k)}, \ k\geq0)$ be the sequence generated by applying a Jacobi-type process to $H$,
\begin{equation}\label{iterH}
H^{(k+1)}=F_k^*H^{(k)}F_k, \quad H^{(0)}=H, \quad k\geq0,
\end{equation}
where $F_k$ are elementary plane matrices acting in the $(i(k),j(k))$ plane, $i(k)<j(k)$.
Suppose that the following assumptions are satisfied.
\begin{itemize}
\item[\textbf{A1}] The pivot strategy is generalized serial
\item[\textbf{A2}] There is a sequence of unitary elementary plane matrices $(U^{(k)}, \ k\geq0)$ such that $$\lim_{k\rightarrow\infty}(F_k-U_k)=0.$$
\item[\textbf{A3}] The diagonal elements of $F_k$ satisfy the condition
    $$\liminf_{k\rightarrow\infty}|f_{i(k)i(k)}^{(k)}|>0.$$
\item[\textbf{A4}] The sequence $(H^{(k)}, \ k\geq0)$ is bounded.
\end{itemize}
Then the following two conditions are equivalent:
\begin{itemize}
\item[(\textbf{i})] $\displaystyle \lim_{k\rightarrow\infty} |h_{i(k)j(k)}^{(k+1)}|=0$,
\item[(\textbf{ii})] $\displaystyle \lim_{k\rightarrow\infty} S(H^{(k)})=0$.
\end{itemize}
\end{Theorem}

\begin{proof}
Obviously, $(ii)$ implies $(i)$. We need to prove that $(i)$ implies $(ii)$.

Transformation matrices $F_k$ from the~\eqref{iterH} can be written as
\[
F_k=U_k+E_k, \quad  k\geq0,
\]
where $U_k$ are unitary matrices from~\textbf{A2}. Assumption~\textbf{A2} implies $E_k = F_k-U_k \rightarrow 0$ as $k\rightarrow\infty$. Therefore, relation~\eqref{iterH} takes the form
\begin{equation}\label{Ek}
H^{(k+1)}=U_k^*H^{(k)}U_k + T_k, \quad k\geq 0.
\end{equation}
To show that
\begin{equation}\label{Ek0}
\lim_{k\rightarrow\infty}T_k=0,
\end{equation}
one has to use assumption~\textbf{A4}. The proof is easy and same as in~\cite[Corollary~5.3]{har+beg-17}.

Applying the vectorization operator $\ve$ to the both sides of~\eqref{Ek} we obtain
\begin{equation}\label{chi}
\chi^{(k+1)}=\mathcal{R}^{(k)}\chi^{(k)} + g^{(k)}, \quad k\geq 0, \quad \lim_{k\rightarrow\infty}g^{(k)}=0.
\end{equation}
Here we have $\chi^{(k)}=\ve(H^{(k)})$, $\chi^{(k+1)}=\ve(H^{(k+1)})$, and
\begin{equation}\label{gk}
g^{(k)}=\ve (T^{(k)}) + (U_k^*H^{(k)}U_k)_{i(k)j(k)}\ve (e_{i(k)}e_{j(k)}^T) + (U_k^*H^{(k)}U_k)_{j(k)i(k)}\ve (e_{j(k)}e_{i(k)}^T),
\end{equation}
where $e_r$ is the $r$th column of $I_n$.
Note that $\ve(e_{i(k)}e_{j(k)}^T)$ (and respectively, $\ve(e_{j(k)}e_{i(k)}^T)$) is the canonical vector in $\mathbf{C}^{2N}$ with only one non-zero element, $1$ at position $\tau(i(k),j(k))$ (resp.\@ $\tau (j(k),i(k))$).
To obtain $\lim_{k\rightarrow\infty}g^{(k)}=0$ we used the relations~\eqref{gk}, \eqref{Ek}, \eqref{Ek0}, and condition (\textbf{i}).

Recall that one sweep (cycle) of the Jacobi-type process has exactly $N=\frac{n(n-1)}{2}$ steps. To count sweeps we use variable $t$.
Set $\chi^{[t]}=\chi^{(tN)}=\ve(H^{(tN)})$, $t\geq 0$. We have
\[
\chi^{[t]}=\mathcal{J}_{\mathcal{O}}^{[t]} \chi^{[t-1]} + g^{[t]}, \quad t\geq 1, \quad \lim_{t\rightarrow\infty}g^{[t]}=0.
\]
Jacobi operator $\mathcal{J}_{\mathcal{O}}^{[t]}=\mathcal{R}^{(tN-1)}\cdots\mathcal{R}^{(t-1)N+1)}\mathcal{R}^{((t-1)N)}$ is defined by the ordering $\mathcal{O}$ of $\mathcal{P}_n$ and by the appropriate Jacobi annihilators. The spectral norm of any Jacobi annihilator of order $n>2$ is equal to $1$. (And it is zero for $n=2$.) Therefore, we immediately obtain
\[
\|g^{[t]} \|_2 \leq \|g^{((t-1)N)}\|_2+\|g^{((t-1)N+1)}\|_2+\cdots+\|g^{(tN-1)}\|_2, \quad t\geq 1.
\]
Hence, $g^{[t]}\rightarrow0$ as $t\rightarrow\infty$ because $\lim_{k\rightarrow\infty}g^{(k)}=0$ from~\eqref{chi}.

Assumption $\textbf{A1}$ implies that $\mathcal{O}\in\mathcal{C}_{sg}^{(n)}$. Suppose that the chain, in the canonical form, connecting $\mathcal{O}$ with some ordering from $\mathcal{C}_{sp}^{(n)}$ contains $d$ shift equivalences. To use Theorem~\ref{tm:sg}, we have to consider $d+1$ successive sweeps.
Using the spectral norm of a Jacobi operator and relation $\lim_{k\rightarrow\infty}g^{[t]}=0$ we obtain
\begin{equation}\label{iterH-eq}
\chi^{[t+d]}=\mathcal{J}_{\mathcal{O}}^{[t+d]}\cdots\mathcal{J}_{\mathcal{O}}^{[t+1]}
\mathcal{J}_{\mathcal{O}}^{[t]}\chi^{[t-1]} + g_{[d+1]}^{[t]}, \quad t\geq 1, \quad \lim_{k\rightarrow\infty}g_{[d+1]}^{[t]}=0,
\end{equation}
where
\[
\|g_{[d+1]}^{[t]}\|_2 \leq \|g^{[t]} \|_2+\|g^{[t+1]} \|_2+\cdots+\|g^{[t+d]} \|_2.
\]

Next we use \textbf{A3}. Recall that $\displaystyle \lim_{k\rightarrow\infty}(F_k-U_k)=0$. Let
\[
2\nu = \liminf_{k\rightarrow\infty}|f_{i(k)i(k)}^{(k)}|>0.
\]
There is $t_0\geq 1$ such that
\[
|u_{i(k)i(k)}^{(k)}|=|u_{j(k)j(k)}^{(k)}|\geq \nu>0, \quad k\geq t_0 N.
\]
Thus,
\[
\mathcal{J}_{\mathcal{O}}^{[t]}\in \Jl_{\!\!\! \mathcal{O}}^{\nu}, \quad t\geq t_0.
\]
Applying Theorem~\ref{tm:sg} to $\mathcal{J}_{\mathcal{O}}^{[t+d]}\cdots\mathcal{J}_{\mathcal{O}}^{[t+1]}
\mathcal{J}_{\mathcal{O}}^{[t]}$ we get
\begin{equation}\label{iterH-jop}
\|\mathcal{J}_{\mathcal{O}}^{[t+d]}\cdots\mathcal{J}_{\mathcal{O}}^{[t+1]}
\mathcal{J}_{\mathcal{O}}^{[t]}\|_2 \leq\zeta_{n,\nu}, \quad 0\leq\zeta_{n,\nu}<1.
\end{equation}
Using the bound from~\eqref{iterH-jop} in~\eqref{iterH-eq} yields
\begin{align*}
\|\chi^{[k+d]}\|_2 & \leq \|\mathcal{J}_{\mathcal{O}}^{[t+d]}\cdots\mathcal{J}_{\mathcal{O}}^{[t+1]}
\mathcal{J}_{\mathcal{O}}^{[t]}\|_2 \|\chi^{[t-1]}\|_2 + \|g_{[d+1]}^{[t]}\|_2 \\
& \leq \zeta_{n,\nu} \|\chi^{[t-1]}\|_2 + \|g_{[d+1]}^{[t]}\|_2.
\end{align*}
Since $0\leq\zeta_{n,\nu}<1$ and $\lim_{t\rightarrow\infty}\|g_{[d+1]}^{[t]}\|_2=0$ we can invoke~\cite[Lemma~1]{har-86} (see also \cite[2.2~Lemma]{har-82}) to obtain $\lim_{t\rightarrow\infty}\chi^{[t]}=0$.

For $(t-1)N< k < tN$, from~\eqref{chi} we obtain
\begin{eqnarray*}
  \|\chi^{(k)}\|_2 & \leq & \|\chi^{((t-1)N)}\|_2+\|g^{((t-1)N)}\|_2+\|g^{((t-1)N+1)}\|_2+\cdots\ +\|g^{(k-1)}\|_2 \\
   & \leq & \|\chi^{[t-1]}\|_2 + (N-1)\max_{(t-1)N<k<tN} \|g^{(k)}\|_2
\end{eqnarray*}
and
$$\|\chi^{(k)}\|_2 \rightarrow0 \quad \text{as} \ t\rightarrow\infty.$$
Therefore, $\lim_{k\rightarrow\infty}\|\chi^{(k)}\|_2=0$. Since $S(H^{(k)})=\|\chi^{(k)}\|_2$, $k\geq 0$,
we have $\lim_{k\rightarrow\infty}S(H^{(k)})=0$.
\end{proof}

The first immediate application of Theorem~\ref{tm:general} is for the complex Jacobi method
for Hermitian matrices (see Corollary~\ref{cor:3.7} above). 
As for the convergence of the diagonal elements, we note that Mascarenhas~\cite{mas-95} has proved that the diagonal elements of the iteration matrix obtained by the Jacobi method always converge. Although his proof is made for the real Jacobi method, it holds for the complex Jacobi method as well.

The second immediate  application of Theorem~\ref{tm:general} is for the complex $J$-Jacobi method which is a natural extension of the real method from \cite{ves-93}. The proof is quite similar to that from \cite[Section~5.1]{har+beg-17}.

\section{An Application to a Jacobi Method for PGEP}\label{sec:pgep}

In this section we show how the new results from previous sections can be used in proving the global convergence of Jacobi methods for the \emph{positive definite  generalized eigenvalue problem} (PGEP)
\[
Ax=\lambda Bx, \quad x\neq0,
\]
where $A,B\in\mathbb{C}^{n\times n}$ are complex Hermitian matrices and $B$ is positive definite.

The block Jacobi methods~\cite{nov+sin-15} seem to be the best methods for solving PGEP with large matrices on contemporary parallel computing machines. At each parallel step they have to solve certain number of PGEPs with smaller matrices, say of dimension $32$--$512$. For such tasks the most appropriate methods are element-wise Jacobi methods. Namely, most of the time, these subproblems have to solve PGEP with almost diagonal matrices. On such matrices element-wise Jacobi methods are very efficient and highly accurate. In this context we say that element-wise Jacobi methods are used as kernel algorithms for the block methods.

There are several Jacobi methods for solving PGEP, such as the complex Falk-Langemeyer (FL) method~\cite{har-CFL}, Hari-Zimmermann (HZ) method~\cite{har-84} and Cholesky-Jacobi (CJ) method~\cite{har-AIP}. Preliminary numerical tests indicate that they all have excellent accuracy properties. Hence, each of these methods is a good candidate for the kernel algorithm of the block methods.
Here we prove the global convergence of the complex CJ method from \cite{har-AIP}. It is the complex extension of the real method which was derived in \cite{har-HZ}. In \cite{har-HZ} the global convergence of the real method is proved. We follow the lines of that proof and, instead of real, we  use  complex Jacobi annihilators and operators from sections 2 and 3. The main role is played by Theorem~\ref{tm:general}.

The complex CJ method has a preliminary step which makes the diagonal elements of the positive definite matrix $B$ ones. Such a step has a preconditioning effect since it most often reduces the condition number of $B$. This is achieved by the congruence transformation with a suitable diagonal matrix $D$,
\[
A^{(0)}=DAD, \quad B^{(0)}=DBD, \quad\ D=\diag(b_{11}^{-\frac{1}{2}},\ldots,b_{nn}^{-\frac{1}{2}}).
\]
Afterwards, in each step $k$ the pivot submatrices of $A^{(k)}$ and  $B^{(k)}$ are simultaneously diagonalized, while maintaining the unit diagonal of $B^{(k)}$. The $k$th step of the CJ method has the form
\begin{equation}\label{PGEPagm}
A^{(k+1)}=Z_k^*A^{(k)}Z_k, \quad B^{(k+1)}=Z_k^*B^{(k)}Z_k, \quad k\geq0,
\end{equation}
where each $Z_k$ is elementary plane matrix that differs from the identity matrix $I_n$ in the $2\times2$ pivot submatrix $\hat{Z}_k$.

We say that the \emph{CJ method for PGEP is convergent} on pair $(A,B)$ if $A^{(k)}\rightarrow\Lambda$ and $B^{(k)}\rightarrow I_n$, where $\Lambda$ is the diagonal matrix of the eigenvalues of $(A,B)$. The method is \emph{globally convergent} if it is convergent on every pair $(A,B)$ with Hermitian $A$ and positive definite $B$.

\subsection{The complex CJ method}

Let us fix $k$ and denote the matrices $Z_k$, $A^{(k)}$ and $B^{(k)}$ from relation~\eqref{PGEPagm} by $Z$, $A$ and $B$, respectively. Their pivot submatrices are denoted by $\hat{Z}$, $\hat{A}$ and $\hat{B}$, respectively. The pivot elements of $A$, $B$ are denoted by  $a_{ij}$, $b_{ij}$, respectively.

The matrix $\hat{Z}$ is obtained as product of the inverse of the Hermitian transpose of the Cholesky factor of $\hat{B}$, denoted by $\hat{C}$, and of the complex Jacobi rotation that diagonalizes $\hat{C}^*\hat{A}\hat{C}$, denoted by $\hat{J}$. Then we have
\begin{equation}\label{CJ-Z}
\hat{Z}=\hat{C}\hat{J}\qquad\mbox{and}\qquad Z=CJ.
\end{equation}
There are two options for choosing $\hat{C}$, one using $LL^*$ factorization and the other using $RR^*$ factorization of $\hat{B}$. The numerical tests from~\cite{har-HZ} indicate that the best accuracy properties of the method are reached when these to approaches are combined in the way that we use
\begin{subequations}
\begin{align}
\hat{B}=\hat{L}\hat{L}^*, & \quad \hat{C}=\hat{L}^{-*}, \quad \text{if} \ a_{ii}^{(k)}\leq a_{jj}^{(k)}, \label{LL*} \\
\hat{B}=\hat{R}\hat{R}^*, & \quad \hat{C}=\hat{R}^{-*}, \quad \text{if} \ a_{ii}^{(k)}> a_{jj}^{(k)}. \label{RR*}
\end{align}
\end{subequations}
Using the first factorization $\hat{B}=\hat{L}\hat{L}^*$, we have
\[
\left[
    \begin{array}{cc}
      1 & b_{ij} \\
      \bar{b}_{ij} & 1 \\
    \end{array}
  \right]=\hat{B}=\hat{L}\hat{L}^*=\left[
    \begin{array}{cc}
      1 & 0 \\
      l_{ji} & l_{jj} \\
    \end{array}
  \right]\left[
    \begin{array}{cc}
      1 & \bar{l}_{ji} \\
      0 & \bar{l}_{jj} \\
    \end{array}
  \right]
\]
implying
\begin{equation}\label{C-L}
\hat{L}=\left[
    \begin{array}{cc}
      1 & 0 \\
      \bar{b}_{ij} & \beta \\
    \end{array}
  \right], \quad
\hat{C}=\hat{L}^{-*}=\frac{1}{\beta}\left[
    \begin{array}{cc}
      \beta & -b_{ij} \\
      0 & 1 \\
    \end{array}
  \right], \quad  \beta=\sqrt{1-|b_{ij}|^2}.
\end{equation}
Matrix $J$ is obtained as the complex Jacobi rotation that diagonalizes $\hat{C}^*\hat{A}\hat{C}$.
Let us assume the following form of $\hat{J}$,
\begin{equation}\label{CJ-J}
\hat{J}=\left[
    \begin{array}{cc}
      \cos\phi & -e^{\imath\alpha}\sin\phi \\
      e^{-\imath\alpha}\sin\phi & \cos\phi \\
    \end{array}
  \right].\end{equation}
A simple calculation using the formulas~\eqref{alpha} and~\eqref{hjacobit} yields (see \cite{har-AIP})
\begin{align*}
\alpha & = \text{arg}(a_{ij}-b_{ij}a_{ii}),  \\
\tan(2\phi) & = \frac{2|a_{ij}-b_{ij}a_{ii}|\beta}
{a_{ii}-a_{jj}+a_{ij}\bar{b}_{ij}+\bar{a}_{ij}b_{ij}-2a_{ii}|b_{ij}|^2}, \quad \phi\in [-\frac{\pi}{4},\frac{\pi}{4}].
\end{align*}
Transformation formulas for the diagonal elements of $A$ are
\begin{align}
a_{ii}' & = a_{ii}+\tan\phi\cdot \frac{|a_{ij}-a_{ii}b_{ij}|}{\beta} , \label{lltaii} \\
a_{jj}' & = a_{jj}-\frac{a_{ij}\bar{b}_{ij}+\bar{a}_{ij}b_{ij}-(a_{ii}+a_{jj})|b_{ij}|^2}{1-|b_{ij}|^2}-\tan\phi\cdot \frac{|a_{ij}-a_{ii}b_{ij}|}{\beta}. \label{lltajj}
\end{align}
If $a_{ii}=a_{jj}$, $a_{ij}=a_{ii}b_{ij}$, the expression for
$\tan (2\phi )$ has the form $0/0$, and then we choose $\phi =0$.
In that case we have $\hat{Z}^*\hat{A}\hat{Z}=a_{ii}I_2$. Hence,
$a_{ii}'$ and $a_{jj}'$ are reduced to $a_{ii}$ and $a_{jj}$, respectively.

On the other hand, if one uses the factorization from~\eqref{RR*}, one obtains (see~\cite{har-AIP})
\begin{equation}\label{C-R}
\hat{C}=\hat{R}^{-*}=
\frac{1}{\beta}\left[\begin{array}{cc}1 & 0 \\ -\bar{b}_{ij} & \beta \end{array}\right],
\end{equation}
and the angles $\alpha$ and $\phi$ in~\eqref{CJ-J} are defined by
\begin{align*}
\alpha & = \text{arg}(a_{ij}-b_{ij}a_{jj}), \\
\tan(2\phi) & = \frac{2|a_{ij}-b_{ij}a_{jj}|\beta}
{a_{ii}-a_{jj}-a_{ij}\bar{b}_{ij}-\bar{a}_{ij}b_{ij}+2a_{jj}|b_{ij}|^2},\quad\ \phi\in[-\frac{\pi}{4},\frac{\pi}{4}].
\end{align*}
The transformation formulas for the diagonal elements of $A$ are
\begin{eqnarray}
a_{ii}' &=& a_{ii}-\frac{a_{ij}\bar{b}_{ij}+\bar{a}_{ij}b_{ij}-(a_{ii}+a_{jj})|b_{ij}|^2}{1-|b_{ij}|^2} +\tan\phi\cdot \frac{|a_{ij}-a_{jj}b_{ij}|}{\beta} , \label{rrtaii}\\[-5pt]
a_{jj}' &=& a_{jj}-\tan\phi\cdot \frac{|a_{ij}-a_{jj}b_{ij}|}{\beta}.\label{rrtajj}
\end{eqnarray}
Again, in the case $a_{ii}=a_{jj}$, $a_{ij}=a_{jj}b_{ij}$, angle $\phi$ is set to zero. Then $a_{ii}'$ and $a_{jj}'$ reduce to $a_{ii}$ and $a_{jj}$, respectively.

Thus, the CJ method is defined by the relations (\ref{CJ-Z})--(\ref{rrtajj}). Our goal is to prove the following theorem.

\begin{Theorem}\label{tm:CJ}
CJ method is globally convergent under the class of generalized serial pivot strategies.
\end{Theorem}

For the proof of Theorem~\ref{tm:CJ} we need  some auxiliary results. First, we want to prove that all matrices $A^{(k)}$, $B^{(k)}$, generated by the method are bounded. That accounts for the assumption \textbf{A4} of Theorem~\ref{tm:general}. Then we want to prove that $b_{i(k)j(k)}^{(k)}$ tends to zero as $k$ increases. After we prove this, the other assumptions of Theorem~\ref{tm:general} will be easy to show.

\begin{Lemma}\label{lema:CJ1}
Let $A$ and $B$ be Hermitian matrices of order $n$ such that $B$ is positive definite.
Let the sequences of matrices $(A^{(k)},\ k\geq 0)$, $(B^{(k)},\ k\geq 0)$ be generated by the CJ
method applied to the pair $(A,B)$ under an arbitrary pivot strategy. Then the assertions (i)--(iii) hold.
\begin{itemize}
\item[(i)] The matrices generated by the method are bounded, that is
\begin{equation}\label{bounds}
\|B^{(k)}\|_2 < n,\qquad \|A^{(k)}\|_2 \leq \mu \|B^{(k)}\|_2 < n\mu ,
\end{equation}
where
$$\mu =\max_{\lambda\in\sigma (A,B)} |\lambda |.$$
Here $\sigma (A,B)$ is the spectrum of  $(A,B)$ and $\mu$ is the spectral radius of the pair $(A,B)$.
\item[(ii)] For the pivot element $b_{i(k)j(k)}^{(k)}$ of $B^{(k)}$ we have
\begin{equation}\label{pivB}
\lim_{k\rightarrow\infty} b_{i(k)j(k)}^{(k)} = 0.
\end{equation}
\item[(iii)] For the transformation matrices $Z_k=C_k J_k$ we have
\[
\lim_{k\rightarrow\infty} \left(Z_k - J_k\right) \rightarrow 0.
\]
\end{itemize}
\end{Lemma}

\begin{proof}
\begin{itemize}
\item[(i)] The proof of~\eqref{bounds} is identical to the proof of~\cite[Lemma~4.1]{har-HZ}. One only has to replace the adjective symmetric by Hermitian.

\item[(ii)] The proof of~\eqref{pivB} follows the lines in the proof of~\cite[Proposition~4.1]{har-HZ}.
Let $B^{(k)}=(b_{rs}^{(k)})$, $k\geq 0$, and
\[
H(B^{(k)})= \frac{\det (B^{(k)})}{b_{11}^{(k)}b_{22}^{(k)}\cdots b_{nn}^{(k)}} = \det (B^{(k)}), \quad k\geq 0.
\]
By the Hadamard inequality we have
\begin{equation}\label{hadamar}
0 < H(B^{(k)}) \leq 1,\quad k\geq 0.
\end{equation}
From~\eqref{CJ-Z} we have $Z_k=C_k J_k$, where each $J_k$ is the complex Jacobi rotation.
Since $\det(J_k)=1$, we get
\[
\det(Z_k) = \det(C_k) = \det(\hat{C}_k) = \frac{1}{\beta_k}, \quad \beta_k =\sqrt{1-|b_{i(k)j(k)}^{(k)}|^2}, \quad k\geq 0.
\]
Therefore,
\[
H(B^{(k+1)})=\det (B^{(k+1)}) = {\det}^2 (C_k) \det(B^{(k)}) = \frac{1}{1-|b_{i(k)j(k)}^{(k)}|^2} H(B^{(k)}), \quad k\geq 0,
\]
that is
\begin{equation}\label{had1}
H(B^{(k)}) = \left(1-|b_{i(k)j(k)}^{(k)}|^2 \right)H(B^{(k+1)}),\quad k\geq 0.
\end{equation}
From the relations~\eqref{had1} and~\eqref{hadamar} we see that
$H(B^{(k)})$ is a nondecreasing sequence of positive real numbers,
bounded above by $1$. It is convergent with limit $\zeta$,
$0<\zeta\leq 1$. By taking the limit on the both sides of the
equation~\eqref{had1}, after cancelation with $\zeta$, we obtain
\[
1 = \lim_{k\rightarrow\infty} \left(1-|b_{i(k)j(k)}^{(k)}|^2 \right)
=1-\lim_{k\rightarrow\infty}|b_{i(k)j(k)}^{(k)}|^2.
\]
This proves~\eqref{pivB}.

\item[(iii)] Recall that each $Z_k$ is the product $C_k J_k$ where $J_k$ is the complex Jacobi rotation.
From the relations~\eqref{LL*}, \eqref{RR*}, \eqref{C-L} and \eqref{C-R} we conclude that
\begin{equation}\label{Ckk}
\hat{C}_k = \left\{ \begin{array}{ll}
\left[ \begin{array}{cc} 1 & -\frac{b_{i(k)j(k)}^{(k)}}{\beta_k} \\ 0 & \frac{1}{\beta_k} \end{array} \right],
& \text{for} \ a_{ii}^{(k)}\leq a_{jj}^{(k)} \\[2ex]
\left[ \begin{array}{cc} \frac{1}{\beta_k} & 0\\ -\frac{\bar{b}_{i(k)j(k)}^{(k)}}{\beta_k} & 1 \end{array} \right],
& \text{for} \ a_{ii}^{(k)}> a_{jj}^{(k)}
 \end{array} \right., \quad \beta_k=\sqrt{1-|b_{i(k)j(k)}^{(k)}|^2}.
\end{equation}
From assertion (ii) we have $b_{i(k)j(k)}^{(k)} \rightarrow 0$ and consequently $\beta_k\rightarrow 1$ as  $k\rightarrow\infty$. Hence, $\hat{C}_k\rightarrow I_2$ i.e.\@ $C_k\rightarrow I_n$, proving the assertion (iii).
\end{itemize}
\end{proof}

Moreover, we need to estimate how close are the diagonal elements of $A^{(k)}$ to the corresponding eigenvalues of the pair $(A^{(k)},B^{(k)})$. Let the eigenvalues of the initial pair $(A,B)$ be in nonincreasing order,
\begin{equation}\label{eigs}
\lambda_1=\cdots =\lambda_{s_1}>\lambda_{s_1+1}=\cdots =\lambda_{s_2}> \ \cdots\ >
\lambda_{s_{p-1}+1}=\cdots =\lambda_{s_p}.
\end{equation}
The case $p=1$ implies $A=\lambda_1B$. Then every nonzero vector is an eigenvector belonging to the only eigenvalue $\lambda_1$.
Let $p>1$.
If in~\eqref{eigs} we set $s_0=0$, we conclude that $n_r=s_r - s_{r-1}$ is the multiplicity of $\lambda_{s_r}$.
Set $\lambda_{s_0}=\lambda_0=\infty$, $\lambda_{s_{p+1}}=-\infty$, and
\[
3\delta_r= \min \{\lambda_{s_{r-1}} - \lambda_{s_{r}},\lambda_{s_{r}} - \lambda_{s_{r+1}}\}, \quad 1\leq r\leq p.
\]
We see that $3\delta_r$ is the absolute gap in the spectrum of $(A,B)$ associated with $\lambda_{s_r}$.
Let
\begin{equation}\label{delta0}
\delta = \min_{1\leq r\leq p} \delta_r, \quad \delta_0=\frac{\delta}{1+\mu^2},
\end{equation}
where $\mu$ is the spectral radius of $(A,B)$. Obviously, $3\delta$ is the minimum absolute gap and $\delta_0<\delta$.

\begin{Lemma}\label{tm:lemma4.3}
Let $A$, $B$ be Hermitian matrices of order $n$ such that $B$ is positive definite with unit diagonal. Let the eigenvalues of $(A,B)$ be ordered as in relation~\eqref{eigs} and $\delta_0$ be as in relation~\eqref{delta0}. If
\[
S(A,B) < \delta_0,
\]
then there is a permutation matrix $P$ such that for the matrix $\tilde{A}=P^TAP =(\tilde{a}_{rt})$ we have
\[
2\sum_{l=1}^n |\tilde{a}_{ll}-\lambda_l|^2 \leq \frac{S^4(A,B)}{\delta_0^2}.
\]
\end{Lemma}

\begin{proof}
The proof is a reformulation of~\cite[Corollary~3.3]{har-86}. Note that the pair $(P^TAP,P^TBP)$ has the same eigenvalues as the pair $(A,B)$ and $P^TBP$ has unit diagonal.
\end{proof}
Now, we are ready to prove  Theorem~\ref{tm:CJ}.

\begin{proof}[Proof of Theorem~\ref{tm:CJ}]
We apply Theorem~\ref{tm:general} to the sequences $(B^{(k)}, \ k\geq 0)$ and $(A^{(k)}, \ k\geq 0)$. In both cases the assumptions \textbf{A1}, \textbf{A2}, \textbf{A4} and condition (\textbf{i}) hold. Indeed, \textbf{A1} is just the selection of the pivot strategy while \textbf{A2} and \textbf{A4} are the assertions (iii) and (i) of Lemma~\ref{lema:CJ1}, respectively. Condition (\textbf{i}) holds because the CJ method diagonalizes the pivot submatrices, thus, $a_{i(k)j(k)}^{(k+1)}=0$ and $b_{i(k)j(k)}^{(k+1)}=0$ holds for all $k\geq 0$.

It remains to prove the assumption \textbf{A3}, that is $\displaystyle \liminf_{k\rightarrow\infty}|z_{i(k)i(k)}^{(k)}|>0$, where we use notation $Z_k =(z_{rs}^{(k)})$. Note that both $Z_k$ and $J_k$ differ from the identity matrix $I_n$ only in the elements of $\hat{Z}_k$ and $\hat{J}_k$. It is sufficient to prove that \textbf{A3} holds for the pivot submatrices.
For $\hat{Z}_k$ we have
\begin{equation}\label{Zdiag}
\hat{Z}_k=\hat{J}_k +(\hat{C}_k-I_2)\hat{J}_k, \quad k\geq 0.
\end{equation}
From relation~\eqref{Ckk} we easily conclude that
\[
\|\hat{C}_k-I_2\|_2^2 =\left(\frac{1}{\beta_k}-1\right)^2+\frac{|b_{i(k),j(k)}^{(k)}|^2}{\beta_k^2}
=\left( \frac{1}{\beta_k^2}+\frac{|b_{i(k)j(k)}^{(k)}|^2}{\beta_k^2(1+\beta_k)^2}\right)
|b_{i(k)j(k)}^{(k)}|^2.
\]
From assertion (ii) of Lemma~\ref{lema:CJ1} we have $b_{i(k),j(k)}^{(k)}\rightarrow 0$, and consequently
$\beta_k\rightarrow 1$ as $k\rightarrow\infty$. Thus, there is a positive integer $k_0$ such that
\[
\|\hat{C}_k-I_2\|_2 \leq 1.1 |b_{i(k)j(k)}^{(k)}|,\qquad k\geq k_0.
\]
Furthermore, from relation~\eqref{Zdiag} we see that the diagonal elements of $\hat{Z}_k$ cannot differ by more than $1.1\cdot|b_{i(k)j(k)}^{(k)}|$ from $\cos\phi_k\geq \sqrt{2}/2$ for $k\geq k_0$.
This implies
\[
\liminf_{k\rightarrow\infty}|z_{i(k)i(k)}^{(k)}|\geq \frac{\sqrt{2}}{2} > 0.
\]
Therefore, both sequences $(S(A^{(k)}), \ k\geq 0)$ and $(S(B^{(k)}) \ k\geq 0)$ converge to zero. Since each $B^{(k)}$ has unit diagonal, it is shown that $B^{(k)}\rightarrow I_n$ as
$k\rightarrow\infty$.

If $\sigma (A,B)$ is a singleton, i.e.\@ if we have $p=1$ in~\eqref{eigs}, the proof is completed.
Namely, if $A=\lambda_1 B$, we have $A^{(k)}=\lambda_1 B^{(k)}$, $k\geq 0$. Recall that in this case the CJ algorithm chooses $\phi_k=0$, $k\geq 0$.

It remains to prove that the diagonal elements of $A^{(k)}$ converge if $p>1$. This comes down to showing that, for a large enough $k$, the diagonal elements of $A^{(k)}$ cannot change their eigenvalue affiliations. We assume that $k_1\geq 1$ is such that
\begin{equation}\label{cond}
S(A^{(k)},B^{(k)}) \leq \delta_0 \quad \text{and} \quad S(B^{(k)}) \leq \frac{1}{2}, \quad k\geq k_1,
\end{equation}
where $\delta_0$ is from~\eqref{delta0}.
Let us fix $k$, $k\geq k_1$, and denote $A^{(k)}$, $B^{(k)}$, $A^{(k+1)}$, $B^{(k+1)}$ by $A=(a_{rs})$, $B=(b_{rs})$, $A=(a_{rs}')$, $B=(b_{rs}')$, respectively. Then from Lemma~\ref{tm:lemma4.3} we have
\begin{equation}\label{discs1}
2\sum_{l=1}^n |a_{ll}-\lambda_{\pi (l)}|^2 \leq \frac{S^4(A,B)}{\delta_0^2}\leq \delta_0^2,
\end{equation}
for some permutation $\pi$ of the set $\{1,2,\ldots\ ,n \}$.
For the two affected diagonal elements of $A$ we have
\[
|a_{ii}-\lambda_{\pi (i)}|\leq \frac{\delta_0}{\sqrt{2}}<\frac{\delta}{\sqrt{2}}, \quad |a_{jj}-\lambda_{\pi (j)}|\leq \frac{\delta_0}{\sqrt{2}}\leq \frac{\delta}{\sqrt{2}},
\]
where $(i,j)$ is the pivot pair. It is sufficient to show that
\begin{equation}\label{toprove}
|a_{ii}'-\lambda_{\pi (i)}| < (3-\frac{1}{\sqrt{2}})\delta,\qquad |a_{jj}'-\lambda_{\pi (j)}|< (3-\frac{1}{\sqrt{2}}) .
\end{equation}
Pair $(A',B')$ satisfies condition~\eqref{cond} and therefore relation~\eqref{discs1} holds for $A'$ with possibly some other permutation. Anyway, the diagonal elements of $A'$ are contained in the disks of radius $\delta_0/\sqrt{2}$ around the eigenvalues, and consequently have to be contained in the larger discs of radius $\delta/\sqrt{2}$. Therefore, condition~\eqref{toprove} precludes the possibility that $a_{ii}$ or $a_{jj}$ changes the disc after the transformation.

If we  prove
\begin{equation}\label{diag_change}
|a_{ii}'-a_{ii}| \leq \sqrt{2}\delta, \quad |a_{jj}'-a_{jj}| \leq \sqrt{2}\delta,
\end{equation}
we shall have
\[
|a_{ii}'-\lambda_{\pi (i)}|\leq |a_{ii}'-a_{ii}|+|a_{ii}-\lambda_{\pi (i)}|\leq \sqrt{2}\delta  +\delta_0/\sqrt{2}=(\sqrt{2}+\frac{1}{\sqrt{2}})\delta < (3-\frac{1}{\sqrt{2}})\delta,
\]
and the same bound will hold for  $|a_{jj}'-\lambda_{\pi (j)}|$. Thus $a_{ii}'$ (resp.\@ $a_{jj}'$) remains in the same disc where $a_{ii}$ (resp.\@ $a_{jj}$) lies.
Let $\xi = \max\{|a_{ii}'-a_{ii}|,|a_{jj}'-a_{jj}|\}$. From the relations~\eqref{lltaii}, \eqref{lltajj}, \eqref{rrtaii}, \eqref{rrtajj}, we have
\[
\xi \leq  \left| \frac{a_{ij}\bar{b}_{ij}+\bar{a}_{ij}b_{ij}-(a_{ii}+a_{jj})|b_{ij}|^2}{1-|b_{ij}|^2}\right|+
\left\{ \begin{array}{ll} \left| \tan\phi\cdot \frac{|a_{ij}-a_{ii}b_{ij}|}{\beta}\right|,
& a_{ii} \leq a_{jj} \\[2ex]
\left|\tan\phi\cdot \frac{|a_{ij}-a_{jj}b_{ij}|}{\beta}\right|,
& a_{ii} > a_{jj}
\end{array} \right. .
\]
From~\eqref{cond} we have
\[
|b_{ij}|\leq \frac{S(B)}{\sqrt{2}}\leq \frac{\sqrt{2}}{4}, \quad \beta \geq \sqrt{\frac{7}{8}}, \quad
\frac{|b_{ij}|}{\beta}\leq \frac{\sqrt{2}}{4}\frac{1}{\sqrt{7/8}}=\frac{1}{\sqrt{7}}<0.378.
\]
Finally we obtain
\begin{eqnarray*}
  \xi &\leq & 2\frac{|a_{ij}|\,|b_{ij}|+\mu |b_{ij}|^2}{\beta^2}+\frac{|a_{ij}|+\mu |b_{ij}|}{\beta} =  \frac{|a_{ij}|+\mu |b_{ij}|}{\beta} \left(1+2\frac{|b_{ij}|}{\beta}\right)\\
  &\leq & \frac{S(A,B)\sqrt{1+\mu^2}}{\sqrt{2}\beta}\left(1+\frac{2}{\sqrt{7}}\right)\leq
  \left(\frac{2}{\sqrt{7}}+\frac{4}{7}\right)\delta < 1.33\delta .
\end{eqnarray*}
Thus, relation~\eqref{diag_change} holds, which completes the proof of the theorem.
\end{proof}

\end{document}